\documentclass[reqno]{amsart}

\usepackage{amssymb}
\usepackage{hyperref}
\usepackage{ifthen}
\usepackage{xargs}
\usepackage[all]{xy}
\usepackage{mathtools}

\hypersetup{colorlinks=true,linkcolor=blue}

\newcommand{\axitem}[1]{\phantomsection \label{ax:#1}}
\newcommand{\axref}[1]{(\hyperref[ax:#1]{#1})}

\newcommand{\newref}[4][]{
\ifthenelse{\equal{#1}{}}{\newtheorem{h#2}[hthm]{#4}}{\newtheorem{h#2}{#4}[#1]}
\expandafter\newcommand\csname r#2\endcsname[1]{#3~\ref{#2:##1}}
\expandafter\newcommand\csname R#2\endcsname[1]{#4~\ref{#2:##1}}
\expandafter\newcommand\csname n#2\endcsname[1]{\ref{#2:##1}}
\newenvironmentx{#2}[2][1=,2=]{
\ifthenelse{\equal{##2}{}}{\begin{h#2}}{\begin{h#2}[##2]}
\ifthenelse{\equal{##1}{}}{}{\label{#2:##1}}
}{\end{h#2}}
}

\newref[section]{thm}{Theorem}{Theorem}
\newref{lem}{Lemma}{Lemma}
\newref{prop}{Proposition}{Proposition}
\newref{cor}{Corollary}{Corollary}
\newref{cond}{Condition}{Condition}

\theoremstyle{definition}
\newref{defn}{Definition}{Definition}
\newref{example}{Example}{Example}

\theoremstyle{remark}
\newref{remark}{Remark}{Remark}

\numberwithin{figure}{section}

\newcommand{\cover}{\triangleleft}
\newcommand{\cat}[1]{\mathbf{#1}}
\renewcommand{\deg}[1]{\mathrm{deg}(#1)}

\DeclarePairedDelimiter\abs{\lvert}{\rvert}

\makeatletter
\let\oldabs\abs
\def\abs{\@ifstar{\oldabs}{\oldabs*}}
\makeatother

\begin{document}

\title{Schemes without Sheaves}

\author{Valery Isaev}

\begin{abstract}
We introduce a new, elementary, and fully constructive approach to schemes based on the notion of a \emph{scheme site}.
A scheme site is a preorder equipped with a presheaf of rings satisfying a simple axiom expressed in terms of localizations.
This framework avoids the use of topological spaces, sheaves, and locally ringed spaces in the foundational definition, while retaining the essential geometric content.

We show that scheme sites form a category.
Scheme sites satisfying an additional saturation condition form a coreflective subcategory, equivalent to the category of schemes.
Consequently, every scheme site presents a scheme, namely the scheme of its saturation.
A key feature of our approach is that the underlying space of a scheme can be recovered from a covering relation determined by the presheaf of rings.
As a consequence, the sheaf condition does not need to be imposed separately, but follows automatically from the algebraic structure.
\end{abstract}

\maketitle

\section{Introduction}

Scheme theory, introduced by Grothendieck in the 1960s \cite{EGA1}, provides a far-reaching generalization of classical algebraic geometry by systematically incorporating nilpotent elements and enabling a flexible functorial viewpoint.
While this framework has proved extraordinarily powerful, it is inherently abstract and technically demanding, and these difficulties arise already at the level of the definition:
a scheme is formulated in terms of several layers of sophisticated concepts, including sheaves, locally ringed spaces, and the spectrum of a ring.

Moreover, even when working with concrete examples such as $\mathrm{Spec}(R)$, one encounters an additional layer of subtlety.
The structure sheaf is defined explicitly only on a basis of open sets,
while its values on arbitrary open subsets are obtained indirectly via the sheaf condition (see \cite[Definition~1.3.4]{EGA1} and \cite[Tag~01HR]{stacks}).
In practice, computations and constructions are carried out almost exclusively on these basic open sets.

Let us recall how this works in the case of $\mathrm{Spec}(R)$.
Consider the standard basis of principal open subsets $D(a) \subseteq \mathrm{Spec}(R)$. The structure sheaf is explicitly given on these opens by
\[ \mathcal{O}(D(a)) = R[a^{-1}], \]
the localization of $R$ at $a$. If $D(b) \subseteq D(a)$, the restriction map
\[ \mathcal{O}(D(a)) \to \mathcal{O}(D(b)) \]
is simply the canonical localization map $R[a^{-1}] \to R[b^{-1}]$.
Thus, on this basis, everything is described in purely algebraic terms using localizations.

On the other hand, for an arbitrary open subset $U \subseteq \mathrm{Spec}(R)$,
the ring $\mathcal{O}(U)$ is not given directly, but rather defined as the set of compatible families of sections over coverings of $U$ by principal opens.
In other words, its definition relies entirely on the sheaf condition, and it does not introduce essentially new algebraic data beyond what is already present on the basic opens.

This observation motivates the main idea of the present paper: the role of the sheaf condition is primarily to ensure that sections over arbitrary open sets are uniquely determined by their values on a generating family.
Consequently, a scheme carries a significant amount of redundant data, since much of the structure is uniquely recoverable from a considerably smaller collection of information.

By eliminating this redundancy, one can not only simplify the underlying data but also streamline several technical aspects of the definition of a scheme.
In this paper, we demonstrate that schemes can be described in a more elementary and convenient way, avoiding much of the traditional formal overhead.

Our approach may be particularly useful for beginners, as it reduces the conceptual barrier to entry, and it is also well-suited for formalization.
In fact, the entire development presented here has been fully formalized in the proof assistant Arend \cite{arend};
the formalization is part of its standard library \cite{arend-lib} (module \texttt{AG.SchemeSite} and its submodules).
Furthermore, our construction is entirely constructive.

The central notion we introduce is that of a \emph{(saturated) scheme site}.
Informally, a scheme site can be thought of as a scheme together with a chosen atlas of affine open subsets.
More precisely, it consists of a preorder equipped with a presheaf of rings satisfying a single simple condition (see \rdefn{scheme-site}).
A scheme site is \emph{saturated} if it satisfies one more condition, which says that it contains enough affine opens (see \rdefn{sat-scheme-site}).
The main objects of interest in this paper are saturated scheme sites; the more general notion of a scheme site is used to construct them through a process called \emph{saturation} (\rdefn{sat}).

From the perspective of scheme sites, the example of $\mathrm{Spec}(R)$ becomes trivial:
it is presented by the scheme site with a single element whose ring of sections is $R$.
Applying the saturation construction to this scheme site, we obtain the following description of $\mathrm{Spec}(R)$ as a saturated scheme site.
Its underlying preorder is $R$, with the relation
\[
b \leq a \quad \text{if and only if} \quad a \mid b.
\]
The structure presheaf $\mathcal{O}$ is defined as
\[
\mathcal{O}_a = R[a^{-1}],
\]
for an element $a \in R$, with restriction maps given by the canonical localizations (see \rexample{spec}).

Every saturated scheme site gives rise to a scheme, and conversely, every scheme can be represented by a saturated scheme site, although this representation is not unique.
We also define morphisms of scheme sites in such a way that morphisms between saturated scheme sites correspond exactly to morphisms of the associated schemes.
As a result, saturated scheme sites form a category equivalent to the category of schemes.
Moreover, this category is a coreflective subcategory of the category of all scheme sites, the coreflector being the saturation.
Thus, scheme sites provide a more elementary and constructive alternative to the classical notion of a scheme.

\subsection{Related work}

The Lean formalization of schemes presented in \cite{schemes-lean} has undergone several iterations.
In that work, the authors observe that it is more convenient to work with an abstract notion of localization of rings rather than with explicit constructions; this insight is equally applicable in our setting.
However, our framework is more elementary, since we do not need to introduce either sheaves or locally ringed spaces at all.
Moreover, our development is entirely constructive.

In constructive approaches to scheme theory, it is common to replace topological spaces with locales \cite{blechschmidt}, formal topologies \cite{schemes-schuster}, or lattices \cite{schemes-zeuner,schemes-coquand}.
The lattice-theoretic approaches are based on the constructive Zariski spectrum, which goes back to \cite{joyal-spectrum}.
We go further by working with arbitrary preorders, making the theory not only constructive, but also considerably simpler.
A different, synthetic approach to constructive algebraic geometry in homotopy type theory is developed in \cite{synthetic-ag}.

\subsection{Structure of the paper}

In Section~\ref{sec:scheme-site}, we introduce the central notion of this paper, namely scheme sites.
We also define saturated scheme sites and the saturation of a scheme site,
and present two fundamental examples: the spectrum of a ring and the projective spectrum of a graded ring.
In Section~\ref{sec:morphisms}, we define morphisms of scheme sites, thereby obtaining the category $\mathbf{SchSite}$ and its full subcategory $\mathbf{SatSchSite}$ of saturated scheme sites.
We also show that the saturation is right adjoint to the inclusion $\mathbf{SatSchSite} \to \mathbf{SchSite}$.
These two sections require only minimal background, namely the notions of presheaves and localizations of rings.

Section~\ref{sec:equivalence} is devoted to proving that $\mathbf{SatSchSite}$ is equivalent to the category of schemes $\mathbf{Sch}$.
This section requires additional notions, in particular locales, sheaves, and locally ringed locales.
In Section~\ref{sec:scheme-site-to-scheme}, we construct the scheme associated to a saturated scheme site.
Section~\ref{sec:affine-opens} is devoted to the characterization of affine opens.
In Section~\ref{sec:functor}, we prove that the construction of Section~\ref{sec:scheme-site-to-scheme} is functorial and show that the resulting functor is fully faithful.
Finally, in Section~\ref{sec:scheme-to-scheme-site}, we complete the proof of the equivalence by constructing a saturated scheme site associated to a scheme.

\section{Scheme sites}
\label{sec:scheme-site}

In this section, we introduce our central notion of scheme sites and construct the fundamental examples of scheme sites,
namely the spectrum of a ring $\mathrm{Spec}(R)$ and the projective spectrum of a graded ring $\mathrm{Proj}(R)$.
The only categorical and algebraic prerequisites required in this development are presheaves on preorders and localizations of rings.
In particular, in this section, we do not assume any prior knowledge of sheaves or locally ringed spaces.
Throughout the paper, all rings are assumed to be commutative with unit, and all ring homomorphisms preserve the unit.

We begin by recalling the necessary definitions.
A \emph{preorder} is a set equipped with a reflexive and transitive relation.
Let $X$ be a preorder.
A \emph{presheaf of rings} $\mathcal{O}^X$ on $X$ consists of the following data:
\begin{itemize}
\item for every $a \in X$, a ring $\mathcal{O}^X_a$,
\item for every $b \leq a$, a ring homomorphism
\[
\mathcal{O}^X_{b,a} : \mathcal{O}^X_a \to \mathcal{O}^X_b,
\]
called the \emph{restriction map},
\end{itemize}
such that
\begin{itemize}
\item for every $a \in P$, we have $\mathcal{O}^X_{a,a} = \mathrm{id}_{\mathcal{O}^X_a}$,
\item for every $c \leq b \leq a$, we have
\[
\mathcal{O}^X_{c,b} \circ \mathcal{O}^X_{b,a} = \mathcal{O}^X_{c,a}.
\]
\end{itemize}

Finally, we recall the notion of localization of a ring.
As observed in~\cite[Section~4.1]{schemes-lean}, it is convenient to work with an abstract characterization of localizations rather than a specific construction such as $R[S^{-1}]$.

\begin{defn}[ring-loc]
Let $f : R \to T$ be a homomorphism of rings.
Let $S$ be a multiplicative submonoid of $R$.
We say that $f$ is a \emph{localization} at $S$ if the following conditions hold:
\begin{enumerate}
\item \label{loc:inv} For every $x \in S$, the element $f(x)$ is invertible.
\item \label{loc:surj} For every $z \in T$, there exist $y \in R$ and $x \in S$ such that $z f(x) = f(y)$.
\item \label{loc:inj} For every $y \in R$, if $f(y) = 0$, then there exists $x \in S$ such that $y x = 0$.
\end{enumerate}
\end{defn}

The last two conditions are equivalent to the usual universal property of localizations:
for every homomorphism of rings $g : R \to T'$ such that $g(x)$ is invertible for every $x \in S$,
there exists a unique homomorphism $\widetilde{g} : T \to T'$ satisfying $\widetilde{g} \circ f = g$.

The localization at $S$ always exists and is unique up to an isomorphism.
The canonical localization $R[S^{-1}]$ is defined as the set of equivalence classes of pairs $(y,x)$, where $y \in R$ and $x \in S$.
Two pairs $(y_1,x_1)$ and $(y_2,x_2)$ are equivalent if there exists $z \in S$ such that $y_1 x_2 z = y_2 x_1 z$.
We denote the equivalence class of a pair $(y,x)$ by $y/x$.
The map $R \to R[S^{-1}]$ is given by $y \mapsto y/1$.
We will work only with localizations at sets of the form $S_x = \{ x^n \mid n \in \mathbb{N} \}$ for some $x \in R$.
The canonical localization $R[S_x^{-1}]$ will be denoted by $R[x^{-1}]$.

We collect several standard facts about localizations that will be used throughout the paper.
Proofs are straightforward and can be found in the literature (see, for example, \cite[Chapter~3]{atiyah-macdonald} and \cite[Tag~00CM]{stacks}, or \cite{constructive-algebra,cacm} for constructive treatments), so we omit them.

\begin{lem}[loc-comp]
Let $f : R \to T$ and $g : T \to U$ be homomorphisms of rings.
If $f$ is a localization at $x \in R$ and $g$ is a localization at $f(y)$ for some $y \in R$, then $g \circ f$ is a localization at $x y$.
\end{lem}

\begin{lem}[loc-factor]
Let $f : R \to T$ and $g : T \to U$ be homomorphisms of rings.
If $f$ is a localization at $S \subseteq R$ and $g \circ f$ is a localization at $x \in R$, then $g$ is a localization at $f(x)$.
\end{lem}

\begin{lem}[loc-transport]
Let $f : R \to T$ be a localization at $x \in R$.
Then $f$ is also a localization at $y \in R$ if and only if there exist $n,k \in \mathbb{N}$ such that
\[ x \mid y^n \quad \text{and} \quad y \mid x^k. \]
\end{lem}

\begin{lem}[loc-ideal]
Let $f : R \to T$ be a localization at $S \subseteq R$.
Let $U$ be a subset of $R$.
Then the ideal of $T$ generated by the set $f(U)$ is the unit ideal if and only if the ideal generated by $U$ contains an element of $S$.
\end{lem}

\subsection{Definition of scheme sites}

We are now ready to introduce the central notion of this paper.

\begin{defn}[scheme-site]
A \emph{scheme site} is a preorder $X$ together with a presheaf of rings $\mathcal{O}^X$ on $X$, called \emph{the structure presheaf}, satisfying the following condition:
\begin{enumerate}
\item[(SC)] \axitem{SC} For every $a',a,b \in X$ such that $a' \leq a$ and $b \leq a$, there exist $x \in \mathcal{O}^X_a$ and $b' \in X$ such that
    $b' \leq a'$, $b' \leq b$, $\mathcal{O}^X_{b,a}$ is a localization at $x$, and $\mathcal{O}^X_{b',a'}$ is a localization at $\mathcal{O}^X_{a',a}(x)$.
\end{enumerate}
We call elements of $X$ \emph{affine opens} of the scheme site.
\end{defn}

The elements of a scheme site $X$ should be thought of as affine open subsets of the corresponding scheme, although not every affine open is required to belong to $X$.
If $b \leq a$, then the affine open corresponding to $b$ is contained within the one corresponding to $a$,
but the converse does not necessarily hold: inclusion of affine opens is only partially reflected by the preorder structure.

Condition \axref{SC} can be visualized as the following diagram, where labels on arrows denote elements with respect to which the corresponding restriction maps are localizations.
\[ \xymatrix{ b' \ar@{-->}[r] \ar@{-->}[d]_{\mathcal{O}^X_{a',a}(x)} & b \ar[d]^x \\
              a' \ar[r]                                              & a
	    } \]

This condition combines two properties.
The first property is that all restriction maps are localizations.
The second property is the existence of intersections: the element $b'$ plays the role of the intersection of $a'$ and $b$.
In fact, we can show that it is equivalent to the combination of these properties:

\begin{prop}
In the context of \rdefn{scheme-site}, condition \axref{SC} is equivalent to the following conditions:
\begin{enumerate}
\item[(SR)] \axitem{SR} For every $a \in X$ and $b \leq a$, the map $\mathcal{O}^X_{b,a} : \mathcal{O}^X_a \to \mathcal{O}^X_b$ is a localization at some $x \in \mathcal{O}^X_a$.
\item[(SS)] \axitem{SS} For every $a',a,b \in X$ such that $a' \leq a$ and $b \leq a$, if $\mathcal{O}^X_{b,a}$ is a localization at some $x \in \mathcal{O}^X_a$,
then there exists $b' \in X$ such that $b' \leq a'$, $b' \leq b$, and $\mathcal{O}^X_{b',a'}$ is a localization at $\mathcal{O}^X_{a',a}(x)$.
\end{enumerate}
\end{prop}
\begin{proof}
Clearly, \axref{SR} and \axref{SS} together imply \axref{SC}.
In the opposite direction, \axref{SR} follows from \axref{SC} if we take $a' = a$.
Finally, \axref{SC} implies \axref{SS} by \rlem{loc-transport}.
\end{proof}

A scheme site is not required to contain all the affine opens of the scheme it presents.
For example, if $a$ is an affine open of a scheme, then every element $x$ of the ring of sections $\mathcal{O}_a$ determines an affine subopen of $a$ with the ring of sections $\mathcal{O}_a[x^{-1}]$.
This subopen is called the \emph{principal open} of $x$ and usually denoted by $D(x)$.
Scheme sites do not have to contain principal opens.
The ones that do contain them play a special role:

\begin{defn}[sat-scheme-site]
A scheme site $X$ is \emph{saturated} if it satisfies the following condition:
\begin{enumerate}
\item[(SE)] \axitem{SE} For every $a \in X$ and $x \in \mathcal{O}^X_a$, there exists $b \leq a$ such that $\mathcal{O}^X_{b,a}$ is a localization at $x$.
\end{enumerate}
\end{defn}

Condition \axref{SE} means that $X$ contains all principal affine opens: the open $b$ plays the role of $D(x)$.
As we will see, only saturated scheme sites give rise to schemes directly.
Nevertheless, every scheme site can be completed to a saturated one (see Section~\ref{sec:saturation}),
and non-saturated scheme sites are often much more economical.

\begin{example}[point]
Let $R$ be a ring.
We denote by $\mathrm{pt}_R$ the scheme site with a single element $\ast$ such that $\mathcal{O}_\ast = R$.
Condition \axref{SC} holds with $x = 1$ and $b' = \ast$.
This scheme site is not saturated unless $R$ is trivial since \axref{SE} implies that every element of $R$ is invertible.
We will see that $\mathrm{pt}_R$ presents the affine scheme $\mathrm{Spec}(R)$.
\end{example}

\subsection{Saturation}
\label{sec:saturation}

We now describe a construction that turns an arbitrary scheme site into a saturated one.
Intuitively, we add all principal opens of affine opens of the scheme site.

\begin{defn}[sat]
Let $X$ be a scheme site.
The \emph{saturation} of $X$ is the scheme site $\mathrm{Sat}(X)$ defined as follows:
\begin{itemize}
\item The elements of $\mathrm{Sat}(X)$ are pairs $(a,x)$, where $a \in X$ and $x \in \mathcal{O}^X_a$.
\item We have $(b,y) \leq (a,x)$ if and only if $b \leq a$ and $\mathcal{O}^X_{b,a}(x)$ divides $y$.
\item The structure presheaf is given by $\mathcal{O}^{\mathrm{Sat}(X)}_{(a,x)} = \mathcal{O}^X_a[x^{-1}]$.
\item If $(b,y) \leq (a,x)$, then the restriction map $\mathcal{O}^{\mathrm{Sat}(X)}_{(b,y),(a,x)} : \mathcal{O}^X_a[x^{-1}] \to \mathcal{O}^X_b[y^{-1}]$
is the unique map extending the composite $\mathcal{O}^X_a \to \mathcal{O}^X_b \to \mathcal{O}^X_b[y^{-1}]$.
It exists since $\mathcal{O}^X_{b,a}(x)$ divides $y$, so its image in $\mathcal{O}^X_b[y^{-1}]$ is invertible.
\end{itemize}
\end{defn}

The relation $\leq$ is clearly reflexive.
It is also transitive: if $(c,z) \leq (b,y) \leq (a,x)$, then $\mathcal{O}^X_{c,a}(x) = \mathcal{O}^X_{c,b}(\mathcal{O}^X_{b,a}(x))$ divides $\mathcal{O}^X_{c,b}(y)$, which divides $z$.
Functoriality of $\mathcal{O}^{\mathrm{Sat}(X)}$ follows from the universal property of localizations.
Note that $\mathcal{O}^{\mathrm{Sat}(X)}_{(a,1)}$ is canonically isomorphic to $\mathcal{O}^X_a$, so $X$ can be thought of as a part of $\mathrm{Sat}(X)$.

To show that $\mathrm{Sat}(X)$ is a saturated scheme site, we need the following lemma that allows us to compute restriction maps of $\mathrm{Sat}(X)$.

\begin{lem}[sat-loc]
Let $X$ be a scheme site, and let $(b,y) \leq (a,x)$ be elements of $\mathrm{Sat}(X)$.
Suppose that $\mathcal{O}^X_{b,a}$ is a localization at $z \in \mathcal{O}^X_a$, and let $w,c \in \mathcal{O}^X_a$ be elements such that $c$ is a power of $z$ and $y \cdot \mathcal{O}^X_{b,a}(c) = \mathcal{O}^X_{b,a}(w)$.
Then $\mathcal{O}^{\mathrm{Sat}(X)}_{(b,y),(a,x)}$ is a localization at $z w / 1$.
\end{lem}
\begin{proof}
Consider the following commutative square:
\[ \xymatrix{ \mathcal{O}^X_a \ar[r]^-{l_x} \ar[d]_{\mathcal{O}^X_{b,a}} & \mathcal{O}^X_a[x^{-1}] \ar[d]^{\mathcal{O}^{\mathrm{Sat}(X)}_{(b,y),(a,x)}} \\
	      \mathcal{O}^X_b \ar[r]_-{l_y} & \mathcal{O}^X_b[y^{-1}].
            } \]

Since $\mathcal{O}^X_{b,a}(c)$ is invertible in $\mathcal{O}^X_b$, the elements $y$ and $\mathcal{O}^X_{b,a}(w)$ differ by a unit.
By \rlem{loc-transport}, $l_y$ is also a localization at $\mathcal{O}^X_{b,a}(w)$.
By \rlem{loc-comp}, the composite $l_y \circ \mathcal{O}^X_{b,a} = \mathcal{O}^{\mathrm{Sat}(X)}_{(b,y),(a,x)} \circ l_x$ is a localization at $z w$.
The claim follows from \rlem{loc-factor}.
\end{proof}

\begin{prop}[sat-saturated]
For every scheme site $X$, its saturation $\mathrm{Sat}(X)$ is a saturated scheme site.
\end{prop}
\begin{proof}
\axref{SE}
Let $(a,x)$ be an element of $\mathrm{Sat}(X)$, and let $y/x^n$ be an element of $\mathcal{O}^X_a[x^{-1}]$.
Then $(a, x y) \leq (a,x)$.
By \rlem{sat-loc} applied to $z = c = 1$ and $w = x y$, the map $\mathcal{O}^{\mathrm{Sat}(X)}_{(a,xy),(a,x)}$ is a localization at $x y / 1$.
Since $x y / 1 = (y/x^n) \cdot (x^{n+1}/1)$ and $x^{n+1}/1$ is invertible, \rlem{loc-transport} implies that this map is also a localization at $y/x^n$.

\axref{SC}
Let $(a',x'), (a,x), (b,y)$ be elements of $\mathrm{Sat}(X)$ such that $(a',x') \leq (a,x)$ and $(b,y) \leq (a,x)$.
By \axref{SC} for $X$, there exist $z \in \mathcal{O}^X_a$ and $b' \in X$ such that $b' \leq a'$, $b' \leq b$,
$\mathcal{O}^X_{b,a}$ is a localization at $z$, and $\mathcal{O}^X_{b',a'}$ is a localization at $\mathcal{O}^X_{a',a}(z)$.
By \rdefn{ring-loc}.\eqref{loc:surj}, there exist $w \in \mathcal{O}^X_a$ and $c = z^n$ such that $y \cdot \mathcal{O}^X_{b,a}(c) = \mathcal{O}^X_{b,a}(w)$.
By \rlem{sat-loc}, $\mathcal{O}^{\mathrm{Sat}(X)}_{(b,y),(a,x)}$ is a localization at $z w / 1$.

Define $y' = \mathcal{O}^X_{b',a'}(x') \cdot \mathcal{O}^X_{b',b}(y)$.
Then $(b',y') \leq (a',x')$ and $(b',y') \leq (b,y)$.
We have
\[ y' \cdot \mathcal{O}^X_{b',a'}(\mathcal{O}^X_{a',a}(c)) = \mathcal{O}^X_{b',a'}(x') \cdot \mathcal{O}^X_{b',b}(y \cdot \mathcal{O}^X_{b,a}(c)) = \mathcal{O}^X_{b',a'}(x' \cdot \mathcal{O}^X_{a',a}(w)). \]
By \rlem{sat-loc} applied to $(b',y') \leq (a',x')$, the map $\mathcal{O}^{\mathrm{Sat}(X)}_{(b',y'),(a',x')}$ is a localization at
\[ \mathcal{O}^X_{a',a}(z) \cdot x' \cdot \mathcal{O}^X_{a',a}(w) / 1 = (x'/1) \cdot \mathcal{O}^{\mathrm{Sat}(X)}_{(a',x'),(a,x)}(z w / 1). \]
Since $x'/1$ is invertible in $\mathcal{O}^X_{a'}[x'^{-1}]$, \rlem{loc-transport} implies that this map is also a localization at $\mathcal{O}^{\mathrm{Sat}(X)}_{(a',x'),(a,x)}(z w / 1)$.
\end{proof}

\begin{example}[spec]
Let $R$ be a ring.
We define \emph{the spectrum of $R$}, denoted $\mathrm{Spec}(R)$, as the saturated scheme site $\mathrm{Sat}(\mathrm{pt}_R)$.
Explicitly, it can be described as follows:
\begin{itemize}
\item The underlying preorder of $\mathrm{Spec}(R)$ is $(R,\leq)$, where $b \leq a$ if and only if $a$ divides $b$.
\item The structure presheaf $\mathcal{O}^R_a$ is defined as $R[a^{-1}]$.
\item If $b \leq a$, then $\mathcal{O}^R_{b,a} : R[a^{-1}] \to R[b^{-1}]$ is defined by the universal property of $R[a^{-1}]$.
\end{itemize}
By \rprop{sat-saturated}, it is a saturated scheme site.
Moreover, \rlem{sat-loc} implies that if $b = a d$, then $\mathcal{O}^R_{b,a}$ is a localization at $b/1$ and also at $d/1$.
\end{example}

In the classical approach, the definition of schemes depends on the construction of $\mathrm{Spec}(R)$.
More precisely, one first introduces locally ringed spaces, then defines $\mathrm{Spec}(R)$ as a particular example of such a space,
and only afterwards defines schemes and verifies that $\mathrm{Spec}(R)$ is indeed a scheme.

In contrast, in our approach the definition of scheme sites is independent of $\mathrm{Spec}(R)$.
As a result, we obtain a more direct and streamlined sequence of definitions.

\subsection{Projective spectrum}

The projective spectrum of a graded ring is another fundamental example of schemes.
To define the corresponding scheme site, we first recall the necessary definitions.

A \emph{graded ring} is a ring $R$ together with a sequence of additive submonoids $R_n \subseteq R$ satisfying the following conditions:
\begin{itemize}
\item Every element of $R$ can be expressed as a sum $x_0 + \ldots + x_k$, where $x_i \in R_i$ for each $i$.
\item If $x_0 + \ldots + x_k = 0$ and $x_i \in R_i$ for each $i$, then $x_i = 0$ for each $i$.
\item For all $x \in R_n$ and $y \in R_m$, we have $xy \in R_{n+m}$.
\end{itemize}
Elements of $R_n$ are called \emph{homogeneous} elements of degree $n$.

If $x$ is a homogeneous element, we will write $\deg{x} \in \mathbb{N}$ for its degree.
Note that the degree of a homogeneous element is not unique, so $\deg{x}$ is not a function of $x$.
Instead, $\deg{x}$ is just a convenient notation for some natural number such that $x \in R_{\deg{x}}$.

Let $S$ be a multiplicative submonoid of a graded ring $R$.
We will write $R[S^{-1}]_0$ for the subring of $R[S^{-1}]$ consisting of elements that can be presented in the form $x/y$ with $x$ and $y$ homogeneous of the same degree.

We list necessary standard lemmas about graded rings without proofs (see, for example, \cite[Tag~00JL]{stacks}).

\begin{lem}
If $a$ is a homogeneous element of degree $n$, then $a^k$ is a homogeneous element of degree $n k$.
\end{lem}

\begin{lem}[homogen-degree]
If $x$ is a homogeneous element of degree $n$ and of degree $m$, then either $x = 0$ or $n = m$.
\end{lem}

\begin{lem}[homogen-div]
If $a$ and $b$ are homogeneous elements such that $a \mid b$, then there exists a homogeneous element $c$ such that $b = a c$.
\end{lem}

\begin{lem}[homogen-loc-inv]
Let $R$ be a graded ring, and let $a \in R$ be a homogeneous element.
If an element $x \in R[a^{-1}]_0$ is invertible in $R[a^{-1}]$, then its inverse belongs to $R[a^{-1}]_0$.
\end{lem}
%

Now, we are ready to define the projective spectrum of a graded ring (for the classical construction, see \cite[Section~2]{EGA2}, \cite[Section~II.2]{hartshorne}, or \cite[Tag~01M3]{stacks}):

\begin{example}[proj]
Let $R$ be a graded ring.
We define \emph{the projective spectrum of $R$}, denoted $\mathrm{Proj}(R)$, as the following saturated scheme site:
\begin{itemize}
\item The underlying preorder of $\mathrm{Proj}(R)$ consists of homogeneous elements of positive degree with $b \leq a$ if and only if $a$ divides $b$.
\item The structure presheaf $\mathcal{O}^{(0)}_a$ is defined as $R[a^{-1}]_0$.
\item If $b \leq a$, then the restriction of $\mathcal{O}^R_{b,a} : R[a^{-1}] \to R[b^{-1}]$ to $R[a^{-1}]_0$ has the image inside $R[b^{-1}]_0$.
Indeed, by \rlem{homogen-div}, there exists a homogeneous element $c$ such that $b = a c$.
If $x/a^n$ is an element of $R[a^{-1}]_0$ such that $\deg{x} = \deg{a^n}$,
then $\mathcal{O}^R_{b,a}(x/a^n) = x c^n / b^n$ belongs to $R[b^{-1}]_0$ since $\deg{x c^n} = \deg{x} + n \deg{c} = n(\deg{a} + \deg{c}) = \deg{b^n}$.

Thus, we can define $\mathcal{O}^{(0)}_{b,a} : R[a^{-1}]_0 \to R[b^{-1}]_0$ as the restriction of $\mathcal{O}^R_{b,a}$.
\end{itemize}
\end{example}

To verify that $\mathrm{Proj}(R)$ is indeed a saturated scheme site, we will need a few simple lemmas:

\begin{lem}[homogen-func-conv]
Let $R$ be a graded ring, and let $b,a \in R$ be homogeneous elements of positive degrees such that $b \leq a$.
Let $x \in R[a^{-1}]$ be an element such that $\mathcal{O}^R_{b,a}(x)$ belongs to $R[b^{-1}]_0$.
Then there exists $y \in R[a^{-1}]_0$ such that $\mathcal{O}^R_{b,a}(x) = \mathcal{O}^{(0)}_{b,a}(y)$.
\end{lem}
\begin{proof}
By \rlem{homogen-div}, there exists a homogeneous element $c$ such that $b = a c$.
Now, if $x = x'/a^n$, then $\mathcal{O}^R_{b,a}(x) = x' c^n / b^n$.
By assumption, there exists a fraction $y'/b^k$ such that $\deg{y'} = \deg{b^k}$ and $x' c^n / b^n = y'/b^k$.
This implies that there exists $m \in \mathbb{N}$ such that $x' c^n b^{k+m} = y' b^{n+m}$.
By \rlem{homogen-div}, there exists a homogeneous $x''$ such that $x'' c^n b^{k+m} = y' b^{n+m}$.
By \rlem{homogen-degree}, either $y' b^{n+m} = 0$ or
\[ \deg{x''} + n \deg{c} + (k + m) \deg{b} = \deg{y'} + (n + m) \deg{b}. \]

In the former case, we have $\mathcal{O}^R_{b,a}(x) = 0$, so we can just take $y = 0$.
In the latter case, define $y$ as $x'' / a^n$.
The equation above implies that $\deg{x''} = \deg{a^n}$, and we have $\mathcal{O}^R_{b,a}(x) = x' c^n / b^n = x'' c^n / b^n = \mathcal{O}^{(0)}_{b,a}(y)$.
\end{proof}

\begin{lem}[proj-loc]
Let $R$ be a graded ring, and let $b,a \in R$ be homogeneous elements of positive degrees such that $b \leq a$.
Let $x \in R[a^{-1}]_0$ be an element such that $\mathcal{O}^R_{b,a}$ is a localization at $x$.
Then $\mathcal{O}^{(0)}_{b,a}$ is also a localization at $x$.
\end{lem}
\begin{proof}
Let us verify conditions of \rdefn{ring-loc}:
\begin{enumerate}
\item Since $\mathcal{O}^R_{b,a}(x)$ is invertible, by \rlem{homogen-loc-inv}, $\mathcal{O}^{(0)}_{b,a}(x)$ is also invertible.
\item This follows from \rlem{homogen-func-conv} and the corresponding property of the localization map $\mathcal{O}^R_{b,a}$.
\item This follows immediately from the corresponding property of the localization map $\mathcal{O}^R_{b,a}$.
\end{enumerate}
\end{proof}

Now, we can prove that $\mathrm{Proj}(R)$ is a saturated scheme site:

\begin{prop}
For every graded ring $R$, the construction of $\mathrm{Prof}(R)$ from \rexample{proj} satisfies the axioms of saturated scheme sites.
\end{prop}
\begin{proof}
We verify conditions \axref{SE} and \axref{SC}.

\axref{SE}
Let $x/a^n$ be an element of $R[a^{-1}]_0$. Define $b$ as $x a$.
As shown in \rexample{spec}, the map $\mathcal{O}^R_{b,a}$ is a localization at $b/1$.
By \rlem{loc-transport}, it is also a localization at $x/a^n$.
By \rlem{proj-loc}, the map $\mathcal{O}^{(0)}_{b,a}$ is also a localization at $x/a^n$.

\axref{SC}
Let $a',a,b$ be homogeneous elements of positive degree such that $a' \leq a$ and $b \leq a$.
Then there exist homogeneous $c,d \in R$ such that $a' = a c$ and $b = a d$.
Define $b'$ as $a c d$.

As shown in \rexample{spec}, the map $\mathcal{O}^R_{b,a}$ is a localization at $b/1$.
By \rlem{loc-transport}, it is also a localization at $b^{\deg{a}}/a^{\deg{b}}$.
By \rlem{proj-loc}, the map $\mathcal{O}^{(0)}_{b,a}$ is also a localization at this element.

Similarly, $\mathcal{O}^R_{b',a'}$ is a localization at $b/1$.
By \rlem{loc-transport}, it is also a localization at $\mathcal{O}^R_{a',a}(b^{\deg{a}}/a^{\deg{b}}) = b^{\deg{a}} c^{\deg{b}}/a'^{\deg{b}}$.
By \rlem{proj-loc}, the map $\mathcal{O}^{(0)}_{b',a'}$ is also a localization at this element.
\end{proof}

\section{Morphisms of scheme sites}
\label{sec:morphisms}

In this section, we define morphisms of scheme sites making them into a category $\cat{SchSite}$.
Saturated scheme sites form a full subcategory $\cat{SatSchSite}$ of this category.
We show that it is coreflective, and later, we will prove that it is equivalent to the category of schemes.

\subsection{Coverings}

In order to define morphisms of scheme sites, we first introduce an auxiliary notion of coverings for affine opens of a scheme site.
This notion will also play an important role in subsequent developments, in particular in the construction of the scheme associated to a scheme site.

\begin{defn}[cover]
Let $(X, \mathcal{O}^X)$ be a scheme site, let $a$ be an affine open of $X$, and let $U$ be a subset of $X$.
We denote by $C(a,U)$ the following set:
\[ \{ x \in \mathcal{O}^X_a \mid \exists c \in U, \exists b \in X, b \leq a, b \leq c, \mathcal{O}^X_{b,a} \text{ is a localization at } x \}. \]
We say that $U$ \emph{covers} $a$, written $a \cover U$, if the ideal of $\mathcal{O}^X_a$ generated by $C(a,U)$ is the unit ideal.
\end{defn}

The intuition behind this definition is simple: if we think of affine opens as open subsets of a scheme,
then an element $a$ is covered by a family $U$ precisely when the corresponding open subset is contained within the union of $U$.

Let us demonstrate how this definition works in the case of affine schemes.
Let $D(a)$ be a principal open subset, and let $U$ be a family of principal opens $D(b)$.
For simplicity, we assume that principal opens $D(b)$ in $U$ are subsets of $D(a)$.
Then $D(a)$ is covered by $U$ if and only if $a$ lies in the radical of the ideal generated by the elements $b \in U$.
This is equivalent to saying that the ideal of the ring $\mathcal{O}_a$ generated by the elements $b/1$ for $b \in U$ is the unit ideal.

Alternatively, instead of the elements $b/1$, one may use elements $x \in \mathcal{O}_a$ such that the map $\mathcal{O}_{b,a}$ is a localization at $x$.
In the case of affine schemes, this leads to an equivalent definition, but this formulation also works with arbitrary schemes.
Therefore, we can use it as the definition of the covering relation for scheme sites.
Since we do not assume that elements of $U$ are less than or equal to $a$,
instead of taking $b$ directly from $U$, we allow $b$ to be an arbitrary lower bound of $a$ and an element of $U$.
Thus, we get the condition from \rdefn{cover}.

\begin{remark}[cover-restrict]
Note that $C(a,U) = C(a, \{ b \in U \mid b \leq a \})$.
Thus, if we have $a \cover U$, we can always replace $U$ with a set consisting of elements that are less than or equal to $a$.
\end{remark}

We now prove that the relation $\cover$ indeed behaves like a covering relation.

\begin{lem}[cover-prop]
The covering relation $\cover$ satisfies the following properties:
\begin{enumerate}
\item \label{cover:refl} If $a \in U$, then $a \cover U$.
\item \label{cover:mono} If $U \subseteq V$ and $a \cover U$, then $a \cover V$.
\item \label{cover:restrict} If $a \cover U$, then $a \cover \{ b \in X \mid b \leq a, \exists c \in U, b \leq c \}$.
\item \label{cover:left} If $a' \leq a$ and $a \cover U$, then $a' \cover U$.
\item \label{cover:trans} If $a \cover V$ and for every $c \in V$ we have $c \cover U$, then $a \cover U$.
\item \label{cover:meet} If $a \cover U$ and $a \cover V$, then $a \cover \{ c \in X \mid \exists u \in U, \exists v \in V, c \leq u, c \leq v \}$.
\end{enumerate}
\end{lem}
\begin{proof}
\

\eqref{cover:refl} This is obvious since $C(a,U)$ contains $1$ whenever $a$ belongs to $U$.

\eqref{cover:mono} This is obvious since $C(a,U) \subseteq C(a,V)$.

\eqref{cover:restrict} This follows from the equality $C(a,U) = C(a, \{ b \in X \mid b \leq a, \exists c \in U, b \leq c \})$, which holds by definition.

\eqref{cover:left} It is enough to show that for every $x \in C(a,U)$, the element $\mathcal{O}^X_{a',a}(x)$ belongs to $C(a',U)$.
Let $x$ be an element of $C(a,U)$.
Then there exist $c \in U$ and $b \in X$ such that $b \leq a$, $b \leq c$, and $\mathcal{O}^X_{b,a}$ is a localization at $x$.
By \axref{SS}, there exists $b' \in X$ such that $b' \leq a'$, $b' \leq b$, and $\mathcal{O}^X_{b',a'}$ is a localization at $\mathcal{O}^X_{a',a}(x)$.
This shows that $\mathcal{O}^X_{a',a}(x)$ belongs to $C(a',U)$.

\eqref{cover:trans} Since $a \cover V$, there exist elements $x_1, \ldots, x_n \in C(a,V)$ such that $1$ is a linear combination of $x_1, \ldots, x_n$.
Then $1$ is also a linear combination of arbitrary powers of these elements.
Therefore, it is enough to show that for every $x \in C(a,V)$, the ideal generated by $C(a,U)$ contains some power of $x$.

Let $x$ be an element of $C(a,V)$.
Then there exist $c \in V$ and $b \in X$ such that $b \leq a$, $b \leq c$, and $\mathcal{O}^X_{b,a}$ is a localization at $x$.
By \rlem{loc-ideal}, the ideal generated by $C(a,U)$ contains a power of $x$ if and only if the ideal generated by $\mathcal{O}^X_{b,a}(C(a,U))$ contains $1$.
Let us denote by $I$ the ideal generated by $\mathcal{O}^X_{b,a}(C(a,U))$.
By \eqref{cover:left}, we know that $b \cover U$.
Thus, it is enough to show that every element of $C(b,U)$ belongs to $I$.

Let $y$ be an element of $C(b,U)$.
Then there exist $c' \in U$ and $b' \in X$ such that $b' \leq b$, $b' \leq c'$, and $\mathcal{O}^X_{b',b}$ is a localization at $y$.
By \rdefn{ring-loc}.\eqref{loc:surj}, we know that $y = \mathcal{O}^X_{b,a}(z) \mathcal{O}^X_{b,a}(x^m)^{-1}$ for some $z \in \mathcal{O}^X_a$ and $m \in \mathbb{N}$.
By \rlem{loc-transport}, the map $\mathcal{O}^X_{b',b}$ is a localization at $\mathcal{O}^X_{b,a}(z)$.
By \rlem{loc-comp}, the map $\mathcal{O}^X_{b',a}$ is a localization at $x z$.
This implies that $x z$ belongs to $C(a,U)$.
It follows that $y$ belongs to $I$.

\eqref{cover:meet} Denote the set in the conclusion by $W$.
By \eqref{cover:restrict}, we have $a \cover \{ u' \mid \exists u \in U, u' \leq u \}$.
By \eqref{cover:trans}, it is enough to show that $u' \cover W$ for every $u' \leq u \in U$.
By \eqref{cover:left}, we have $u' \cover V$.
By \eqref{cover:restrict}, we have $u' \cover \{ v' \mid v' \leq u', \exists v \in V, v' \leq v \}$.
Every element of the latter set belongs to $W$, which completes the proof by \eqref{cover:mono}.
\end{proof}

The following lemma proves that the structure presheaf is always separated with respect to the covering relation.
It can be shown that it is even a sheaf, but we do not need this fact yet.

\begin{lem}[str-sep]
Suppose $a \cover U$, and let $x,y \in \mathcal{O}^X_a$ be elements such that for all $c \in U$ and $b \in X$ satisfying $b \leq a$ and $b \leq c$, we have $\mathcal{O}^X_{b,a}(x) = \mathcal{O}^X_{b,a}(y)$.
Then $x = y$.
\end{lem}
\begin{proof}
Since $a \cover U$, the unit ideal of $\mathcal{O}^X_a$ is generated by some elements $z_1, \ldots, z_n \in C(a,U)$.
This means that for each $i$, we have $c_i \in U$ and $b_i \in X$ such that $b_i \leq a$, $b_i \leq c_i$, and $\mathcal{O}^X_{b_i,a}$ is a localization at $z_i$.
By assumption, we know that $\mathcal{O}^X_{b_i,a}(x) = \mathcal{O}^X_{b_i,a}(y)$.
By \rdefn{ring-loc}.\eqref{loc:inj}, there exists $k_i \in \mathbb{N}$ such that $(x - y) z_i^{k_i} = 0$.
Since $z_1^{k_1}, \ldots, z_n^{k_n}$ generate the unit ideal, we get that $x - y = 0$.
\end{proof}

We will also need the following simple criterion for coverings.

\begin{lem}[inv-cover]
Let $b,c \leq a$ be elements of a scheme site $X$, and let $x \in \mathcal{O}^X_a$ be an element such that $\mathcal{O}^X_{c,a}$ is a localization at $x$.
If $\mathcal{O}^X_{b,a}(x)$ is invertible, then $b \cover \{ c \}$.
\end{lem}
\begin{proof}
By \axref{SS}, there exists $d \in X$ such that $d \leq b$, $d \leq c$, and $\mathcal{O}^X_{d,b}$ is a localization at $\mathcal{O}^X_{b,a}(x)$.
Since this element is invertible, $\mathcal{O}^X_{d,b}$ is also a localization at $1$.
It follows that $1$ belongs to $C(b, \{ c \})$ and that $b \cover \{ c \}$.
\end{proof}

\subsection{Premorphisms of scheme sites}

We begin by introducing the notion of a premorphism of scheme sites.
Later, morphisms of scheme sites will be defined either as equivalence classes of premorphisms or as premorphisms satisfying a certain completeness condition.
Thus, premorphisms can be viewed as a convenient way of representing morphisms.

Let $X$ and $Y$ be scheme sites.
Naively, one might try to define a (pre)morphism between $X$ and $Y$ as a pair $(f^*,f^\sharp)$,
where $f^* : Y \to X$ and $f^\sharp$ is a natural transformation between $\mathcal{O}^Y$ and $\mathcal{O}^X \circ f^*$.
This definition is too restrictive since $f^*$ always maps affine opens to affine opens, but arbitrary morphisms of schemes do not have to satisfy this property.
Therefore, instead of taking a function $f^* : Y \to X$, we require only a relation $\leq_f$ on $X \times Y$.
Intuitively, $b \leq_f a$ holds if and only if the image of $b$ is contained within $a$.

\begin{defn}[ssm]
Let $X$ and $Y$ be scheme sites.
A \emph{premorphism} $f$ between $X$ and $Y$ is a pair $(\leq_f,f^\sharp)$, where $\leq_f$ is a relation on $X \times Y$ and $f^\sharp_{b,a} : \mathcal{O}^Y_a \to \mathcal{O}^X_b$ is a ring homomorphism defined for all $b,a$ satisfying $b \leq_f a$.
This pair is required to satisfy the following conditions:
\begin{enumerate}
\item \label{ssm:comp} For all $b,b' \in X$ and $a,a' \in Y$, if $b' \leq b$, $b \leq_f a$, and $a \leq a'$, then $b' \leq_f a'$.
\item \label{ssm:nat} For all $b,b' \in X$, $a,a' \in Y$, and $x \in \mathcal{O}^Y_{a'}$, if $b' \leq b$, $b \leq_f a$, and $a \leq a'$, then $f^\sharp_{b',a}(\mathcal{O}^Y_{a,a'}(x)) = \mathcal{O}^X_{b',b}(f^\sharp_{b,a'}(x))$.
\item \label{ssm:top} For every $b \in X$, we have
\[ b \cover \{ b' \in X \mid \exists a \in Y, b' \leq_f a \}. \]
\item \label{ssm:meet} For all $b \in X$ and $a_1,a_2 \in Y$, if $b \leq_f a_1$ and $b \leq_f a_2$, then
\[ b \cover \{ b' \in X \mid \exists a \in Y, a \leq a_1, a \leq a_2, b' \leq_f a \}. \]
\item \label{ssm:loc} For all $b \in X$ and $a,a' \in Y$, if $a' \leq a$ and $b \leq_f a$, then $b \cover L_f(a,a')$, where $L_f(a,a')$ is the set of elements $b_0 \in X$ such that
$b_0 \leq_f a$ and there exist $b' \in X$ and $x \in \mathcal{O}^Y_a$ satisfying $b' \leq b_0$, $b' \leq_f a'$, $\mathcal{O}^Y_{a',a}$ is a localization at $x$, and $\mathcal{O}^X_{b',b_0}$ is a localization at $f^\sharp_{b_0,a}(x)$.
\end{enumerate}
\end{defn}

Condition~\eqref{ssm:comp} expresses a natural compatibility of the relation $\leq_f$ with the preorder structures on $X$ and $Y$.
Condition~\eqref{ssm:nat} expresses the naturality of $f^\sharp$.

Conditions~\eqref{ssm:top} and~\eqref{ssm:meet} require more discussion.
Informally, they express that the inverse image operation associated to $f$ preserves finite meets.
Since the underlying preorders of scheme sites are not assumed to admit finite meets, this property cannot be formulated directly.
Instead, we impose covering conditions that serve as a substitute.

These conditions are closely related to the notion of flatness for functors between sites (see, for example, \cite[Definition~4.1]{shulman-flat}).
Indeed, one can associate to every scheme site a site in the usual sense, and a premorphism of scheme sites induces a functor between the corresponding sites.
Conditions~\eqref{ssm:top} and~\eqref{ssm:meet} are precisely those needed to ensure that this functor is flat.
In this paper, however, we do not use this approach.

Finally, condition~\eqref{ssm:loc} is formally similar to condition~\axref{SC} in the definition of scheme sites.
Elements $b_0 \in L_f(a,a')$ are those for which the following diagram can be completed, where the horizontal arrows correspond to the relation $\leq_f$ rather than $\leq$:
\[ \xymatrix{ b' \ar@{-->}[r] \ar@{-->}[d]_{f^\sharp_{b_0,a}(x)} & a' \ar[d]^x \\
              b_0 \ar[r]                                       & a
            } \]
From a more geometric perspective, this condition corresponds to the locality condition for morphisms of locally ringed locales.
Note that we do not require that such a diagram exists for $b$ itself, but only that $b$ is covered by elements for which it exists.
The stronger condition would be too restrictive if the domain $X$ is not saturated:
an element $b' \leq b$ such that $\mathcal{O}^X_{b',b}$ is a localization at a given element $f^\sharp_{b,a}(x)$ might not exist at all.

The element $x$ in condition~\eqref{ssm:loc} can be chosen arbitrarily, and we may also assume that the elements of the covering are less than or equal to $b$:

\begin{lem}[loc-forall]
Let $f : X \to Y$ be a premorphism of scheme sites.
Let $b \in X$ and $a,a' \in Y$ be such that $a' \leq a$ and $b \leq_f a$, and let $x \in \mathcal{O}^Y_a$ be an element such that $\mathcal{O}^Y_{a',a}$ is a localization at $x$.
Then $b$ is covered by the set of elements $e \leq b$ for which there exists $e' \in X$ such that $e' \leq e$, $e' \leq_f a'$, and $\mathcal{O}^X_{e',e}$ is a localization at $f^\sharp_{e,a}(x)$.
\end{lem}
\begin{proof}
By condition~\eqref{ssm:loc} and \rlem{cover-prop}.\eqref{cover:restrict}, $b$ is covered by the set of elements $e \leq b$ such that $e \leq b_0$ for some $b_0 \in L_f(a,a')$.
Let $e$ be such an element, and let $b'$ and $y$ be the elements corresponding to $b_0$.
Since $\mathcal{O}^Y_{a',a}$ is a localization at both $x$ and $y$, by \rlem{loc-transport}, $x$ divides a power of $y$ and $y$ divides a power of $x$.
The homomorphism $f^\sharp_{b_0,a}$ preserves these divisibility relations.
Thus, \rlem{loc-transport} implies that $\mathcal{O}^X_{b',b_0}$ is also a localization at $f^\sharp_{b_0,a}(x)$.
By \axref{SS}, there exists $e' \in X$ such that $e' \leq e$, $e' \leq b'$, and $\mathcal{O}^X_{e',e}$ is a localization at $\mathcal{O}^X_{e,b_0}(f^\sharp_{b_0,a}(x)) = f^\sharp_{e,a}(x)$.
Since $e' \leq b' \leq_f a'$, we have $e' \leq_f a'$.
\end{proof}

The following lemma shows that premorphisms of scheme sites preserve the covering relation in a certain sense:

\begin{lem}[map-cover]
Let $f : X \to Y$ be a premorphism of scheme sites.
Suppose $b \leq_f a \cover U$.
Then we have
\[ b \cover \{ b' \in X \mid \exists a' \in U, b' \leq_f a' \}. \]
\end{lem}
\begin{proof}
Denote the set in the conclusion by $T$.
Since $a \cover U$, there exist $y_1, \ldots, y_n \in C(a,U)$ such that $1$ is a linear combination of $y_1, \ldots, y_n$.
For each $i$, there exist $c_i \in U$ and $a_i \in Y$ such that $a_i \leq a$, $a_i \leq c_i$, and $\mathcal{O}^Y_{a_i,a}$ is a localization at $y_i$.
By \rlem{loc-forall}, $b$ is covered by the set $V_i$ of elements $e \leq b$ for which there exists $e' \leq e$ such that $e' \leq_f a_i$ and $\mathcal{O}^X_{e',e}$ is a localization at $f^\sharp_{e,a}(y_i)$.
Applying \rlem{cover-prop}.\eqref{cover:meet} several times (together with $b \cover \{ b \}$), we see that $b$ is covered by the set $W$ of elements $c \leq b$ such that for every $i$, there exists $e_i \in V_i$ satisfying $c \leq e_i$.
By \rlem{cover-prop}.\eqref{cover:trans}, it is enough to show that $c \cover T$ for every $c \in W$.

Fix $c \in W$. Then $c \leq_f a$.
For each $i$, choose $e_i \in V_i$ such that $c \leq e_i$, and let $e'_i$ be the corresponding element.
By \axref{SS}, there exists $d_i \in X$ such that $d_i \leq c$, $d_i \leq e'_i$, and $\mathcal{O}^X_{d_i,c}$ is a localization at $\mathcal{O}^X_{c,e_i}(f^\sharp_{e_i,a}(y_i)) = f^\sharp_{c,a}(y_i)$.
Since $d_i \leq e'_i \leq_f a_i \leq c_i$, we have $d_i \in T$.
Thus, $f^\sharp_{c,a}(y_i) \in C(c,T)$.
This completes the proof since $1$ is a linear combination of $f^\sharp_{c,a}(y_1), \ldots, f^\sharp_{c,a}(y_n)$.
\end{proof}

We will also need the following consequence of condition~\eqref{ssm:loc}:

\begin{lem}[restrict]
Let $f : X \to Y$ be a premorphism of scheme sites.
Let $b \in X$ and $a',a \in Y$ be such that $a' \leq a$ and $b \leq_f a$, and let $x \in \mathcal{O}^Y_a$ be an element such that $\mathcal{O}^Y_{a',a}$ is a localization at $x$.
If $f^\sharp_{b,a}(x)$ is invertible, then $b \cover \{ b' \mid b' \leq_f a' \}$.
\end{lem}
\begin{proof}
By \rlem{loc-forall}, $b$ is covered by the set of elements $e \leq b$ for which there exists $e' \leq e$ such that $e' \leq_f a'$ and $\mathcal{O}^X_{e',e}$ is a localization at $f^\sharp_{e,a}(x)$.
Since $f^\sharp_{e,a}(x) = \mathcal{O}^X_{e,b}(f^\sharp_{b,a}(x))$ is invertible, \rlem{inv-cover} implies that $e \cover \{ e' \}$.
Now, the claim follows from \rlem{cover-prop}.\eqref{cover:trans} and \eqref{cover:mono}.
\end{proof}

We now define the identity premorphism and the composition of premorphisms, and show that scheme sites together with premorphisms form a category.

\begin{defn}
Let $X$ be a scheme site. The \emph{identity premorphism} $\mathrm{id}_X$ is defined by
\[
b \leq_{\mathrm{id}_X} a \;\;\text{if and only if}\;\; b \leq a,
\]
and
\[
(\mathrm{id}_X)^\sharp_{b,a} = \mathcal{O}^X_{b,a}.
\]
\end{defn}

It is straightforward to see that $\mathrm{id}_X$ satisfies the axioms of premorphisms.

\begin{defn}
Let $f : X \to Y$ and $g : Y \to Z$ be premorphisms. Their \emph{composition} $g \circ f : X \to Z$ is defined as follows:
\begin{itemize}
\item For $c \in X$ and $a \in Z$, we set $c \leq_{g \circ f} a$ if there exists $b \in Y$ such that
\[
c \leq_f b \quad \text{and} \quad b \leq_g a.
\]
\item For such $c,b,a$, define
\[
(g \circ f)^\sharp_{c,a} = f^\sharp_{c,b} \circ g^\sharp_{b,a}.
\]
\end{itemize}
\end{defn}

\begin{lem}[comp-well]
The composition is well-defined, that is, the map $(g \circ f)^\sharp_{c,a}$ does not depend on the choice of $b$.
\end{lem}
\begin{proof}
Let $b_1,b_2 \in Y$ be a pair of elements such that $c \leq_f b_1 \leq_g a$ and $c \leq_f b_2 \leq_g a$.
By \rdefn{ssm}.\eqref{ssm:meet}, we have
\[ c \cover \{ c' \mid \exists b, b \leq b_1, b \leq b_2, c' \leq_f b \}. \]
By \rlem{str-sep}, it is enough to show that $\mathcal{O}^X_{c',c} \circ f^\sharp_{c,b_1} \circ g^\sharp_{b_1,a} = \mathcal{O}^X_{c',c} \circ f^\sharp_{c,b_2} \circ g^\sharp_{b_2,a}$
for every $c' \leq c$ such that there exists $b$ satisfying $b \leq b_1$, $b \leq b_2$, and $c' \leq_f b$.

Consider the following diagram:
\[ \xymatrix{                                                                       & \mathcal{O}^Y_{b_1} \ar[rr]^{\mathcal{O}^Y_{b,b_1}} \ar[dr]^{f^\sharp_{c,b_1}}  & & \mathcal{O}^Y_b \ar[rd]^{f^\sharp_{c',b}} \\
              \mathcal{O}^Z_a \ar[ur]^{g^\sharp_{b_1,a}} \ar[dr]_{g^\sharp_{b_2,a}} &                                                                                 & \mathcal{O}^X_c \ar[rr]^{\mathcal{O}^X_{c',c}} & & \mathcal{O}^X_{c'} \\
                                                                                    & \mathcal{O}^Y_{b_2} \ar[rr]_{\mathcal{O}^Y_{b,b_2}}  \ar[ur]_{f^\sharp_{c,b_2}} & & \mathcal{O}^Y_b \ar[ur]_{f^\sharp_{c',b}}  
            } \]
The outer hexagon commutes by naturality of $g^\sharp$.
The parallelograms commute by naturality of $f^\sharp$.
This implies the required equality.
\end{proof}

\begin{lem}
The composition satisfies the axioms of premorphisms.
\end{lem}
\begin{proof}
Let $f : X \to Y$ and $g : Y \to Z$ be premorphisms of scheme sites.
Conditions \eqref{ssm:comp} and \eqref{ssm:nat} easily follow from the corresponding properties for $f$ and $g$.
Let us verify condition~\eqref{ssm:top}.

We need to show that $c \cover \{ c'' \mid \exists a \in Z, c'' \leq_{g \circ f} a \}$.
By condition~\eqref{ssm:top} for $f$, we have $c \cover \{ c' \mid \exists b \in Y, c' \leq_f b \}$.
By \rlem{cover-prop}.\eqref{cover:trans}, it is enough to show that for all $b \in Y$ and $c' \leq_f b$, we have $c' \cover \{ c'' \mid \exists a \in Z, c'' \leq_{g \circ f} a \}$.
By condition~\eqref{ssm:top} for $g$, we have $b \cover \{ b' \mid \exists a \in Z, b' \leq_g a \}$.
By \rlem{map-cover}, we get the required property:
\[ c' \cover \{ c'' \mid \exists b' \in Y, \exists a \in Z, c'' \leq_f b', b' \leq_g a \} = \{ c'' \mid \exists a \in Z, c'' \leq_{g \circ f} a \}. \]

Condition~\eqref{ssm:meet} is verified similarly.

Finally, let us verify condition~\eqref{ssm:loc}.
Let $c \leq_f b \leq_g a$ and $a' \leq a$.
By condition~\eqref{ssm:loc} for $g$, we have $b \cover L_g(a,a')$.
By \rlem{map-cover}, $c$ is covered by the set of elements $d$ such that $d \leq_f t$ for some $t \in L_g(a,a')$.
By \rlem{cover-prop}.\eqref{cover:trans}, it is enough to show that every such $d$ is covered by $L_{g \circ f}(a,a')$.
So, let $d \leq_f t$, where $t \leq_g a$, $t' \leq t$, $t' \leq_g a'$, $\mathcal{O}^Z_{a',a}$ is a localization at $x$, and $\mathcal{O}^Y_{t',t}$ is a localization at $g^\sharp_{t,a}(x)$.
By \rlem{loc-forall}, $d$ is covered by the set of elements $e \leq d$ for which there exists $s \leq e$ such that $s \leq_f t'$ and $\mathcal{O}^X_{s,e}$ is a localization at $f^\sharp_{e,t}(g^\sharp_{t,a}(x)) = (g \circ f)^\sharp_{e,a}(x)$.
Since $e \leq_{g \circ f} a$ and $s \leq_{g \circ f} a'$, every such $e$ belongs to $L_{g \circ f}(a,a')$.
\end{proof}

It is straightforward to verify that the identity premorphisms and the composition defined above satisfy the usual associativity and unit laws.
Thus, scheme sites together with premorphisms form a category.

\subsection{Morphisms of scheme sites}

As we will see later, every premorphism of scheme sites induces a morphism between the corresponding schemes.
However, different premorphisms may give rise to the same morphism.
The reason for this is that condition~\eqref{ssm:comp} in the definition of a premorphism is too weak.

More precisely, this condition ensures that if $b \leq b'$ and $b' \leq_f a$, then $b \leq_f a$.
However, it may happen that $b$ is covered by elements $b'$ such that $b' \leq_f a$, without $b \leq_f a$ itself holding.
From a geometric perspective, this situation should still imply that $b \leq_f a$.
One way to address this issue would be to define morphisms of scheme sites as premorphisms satisfying an additional completeness condition ensuring that $\leq_f$ is closed under coverings.
In this paper, however, we take a simpler approach.

We define morphisms of scheme sites as equivalence classes of premorphisms.
For two premorphisms $f$ and $g$ to be equivalent, we do not require the relations $\leq_f$ and $\leq_g$ themselves to coincide, but only their closures under coverings.
In this way, premorphisms can be viewed as presentations of actual morphisms, similarly to how scheme sites present schemes.

\begin{defn}
Let $f,g : X \to Y$ be premorphisms of scheme sites.
We say that $f$ and $g$ are \emph{equivalent}, written $f \sim g$, if the following conditions hold for all $b \in X$ and $a \in Y$:
\begin{align*}
b \leq_f a & \Rightarrow b \cover \{ b' \mid b' \leq_g a \}, \\
b \leq_g a & \Rightarrow b \cover \{ b' \mid b' \leq_f a \}, \\
b \leq_f a, b \leq_g a & \Rightarrow f^\sharp_{b,a} = g^\sharp_{b,a}.
\end{align*}
A \emph{morphism} of scheme sites is an equivalence class of premorphisms under the relation $\sim$.
\end{defn}

\begin{lem}
The relation $\sim$ on premorphisms is an equivalence relation.
\end{lem}
\begin{proof}
Reflexivity and symmetry are clear.
Let us verify transitivity.
Suppose $f \sim g$ and $g \sim h$.
The first two implications follow from \rlem{cover-prop}.\eqref{cover:trans}.
The only non-trivial part is to prove the implication
\[ b \leq_f a, b \leq_h a \Rightarrow f^\sharp_{b,a} = h^\sharp_{b,a}. \]

Assume that $b \leq_f a$ and $b \leq_h a$.
By \rlem{cover-prop}.\eqref{cover:restrict} and condition~\eqref{ssm:comp}, we have $b \cover \{ b' \leq b \mid b' \leq_g a \}$.
By \rlem{str-sep}, it is enough to prove $\mathcal{O}^X_{b',b} \circ f^\sharp_{b,a} = \mathcal{O}^X_{b',b} \circ h^\sharp_{b,a}$ for every $b'$ such that $b' \leq b$ and $b' \leq_g a$.
So, fix such an element $b'$.
Then we have the required equality by the naturality of $f^\sharp$ and $h^\sharp$:
\[ \mathcal{O}^X_{b',b} \circ f^\sharp_{b,a} = f^\sharp_{b',a} = g^\sharp_{b',a} = h^\sharp_{b',a} = \mathcal{O}^X_{b',b} \circ h^\sharp_{b,a}. \]
\end{proof}

We now show that the composition is compatible with this equivalence relation.

\begin{lem}
The composition of premorphisms respects the equivalence relation.
\end{lem}
\begin{proof}
First, we verify that the composition respects $\sim$ on the left.
The first two implications easily follow from \rlem{map-cover}.
The proof of the third implication is identical to the proof of \rlem{comp-well}.


Now, we verify that the composition respects $\sim$ on the right.
Let $f,g : X \to Y$ be equivalent premorphisms of scheme sites, and let $h : Y \to Z$ be another premorphism.
The first two implications are straightforward.
We prove the last one.

Assume $c \leq_{h \circ f} a$ and $c \leq_{h \circ g} a$.
Then there exists $b \in Y$ such that $c \leq_f b$ and $b \leq_h a$.
Since $f \sim g$, we have $c \cover \{ c' \leq c \mid c' \leq_g b \}$.
By \rlem{str-sep}, it is enough to show that for every $c' \leq c$ such that $c' \leq_g b$, we have $f^\sharp_{c',b} \circ h^\sharp_{b,a} = g^\sharp_{c',b} \circ h^\sharp_{b,a}$,
but this follows from the equality $f^\sharp_{c',b} = g^\sharp_{c',b}$.
\end{proof}

It follows that scheme sites together with their morphisms form a category, which we denote by $\cat{SchSite}$.
We denote by $\cat{SatSchSite}$ the full subcategory of $\cat{SchSite}$ consisting of saturated scheme sites.
We will usually denote a morphism by one of its representatives.

\subsection{Saturation as a coreflection}

In this subsection, we show that $\cat{SatSchSite}$ is a coreflective subcategory of $\cat{SchSite}$ and that the coreflector is given by the saturation.

Let $X$ be a scheme site.
We define a premorphism $\varepsilon_X : \mathrm{Sat}(X) \to X$ as follows:
\begin{align*}
(b,y) \leq_{\varepsilon_X} a & \text{ if and only if } b \leq a, \\
(\varepsilon_X)^\sharp_{(b,y),a} & = l_y \circ \mathcal{O}^X_{b,a},
\end{align*}
where $l_y : \mathcal{O}^X_b \to \mathcal{O}^X_b[y^{-1}]$ is the canonical map.
Conditions~\eqref{ssm:comp}, \eqref{ssm:nat}, \eqref{ssm:top}, and \eqref{ssm:meet} of \rdefn{ssm} are obvious.
To verify condition~\eqref{ssm:loc}, let $(b,y) \leq_{\varepsilon_X} a$ and $a' \leq a$.
By \axref{SC}, there exist $b' \in X$ and $x \in \mathcal{O}^X_a$ such that $b' \leq b$, $b' \leq a'$, $\mathcal{O}^X_{a',a}$ is a localization at $x$, and $\mathcal{O}^X_{b',b}$ is a localization at $\mathcal{O}^X_{b,a}(x)$.
Then $(b', \mathcal{O}^X_{b',b}(y)) \leq (b,y)$, and \rlem{sat-loc} (applied to $z = \mathcal{O}^X_{b,a}(x)$, $c = 1$, and $w = y$) implies that
$\mathcal{O}^{\mathrm{Sat}(X)}_{(b',\mathcal{O}^X_{b',b}(y)),(b,y)}$ is a localization at $\mathcal{O}^X_{b,a}(x) y / 1$.
Since $y/1$ is invertible, by \rlem{loc-transport}, it is also a localization at $l_y(\mathcal{O}^X_{b,a}(x)) = (\varepsilon_X)^\sharp_{(b,y),a}(x)$.
Thus, $(b,y) \in L_{\varepsilon_X}(a,a')$.

Now, let $X$ be a saturated scheme site, let $Y$ be a scheme site, and let $f : X \to Y$ be a premorphism.
We define a premorphism $\widehat{f} : X \to \mathrm{Sat}(Y)$ as follows:
\begin{itemize}
\item $b \leq_{\widehat{f}} (a,y)$ if and only if $b \leq_f a$ and $f^\sharp_{b,a}(y)$ is invertible;
\item $\widehat{f}^\sharp_{b,(a,y)} : \mathcal{O}^Y_a[y^{-1}] \to \mathcal{O}^X_b$ is the unique homomorphism extending $f^\sharp_{b,a}$.
\end{itemize}

\begin{lem}[lift]
The pair $\widehat{f}$ is a premorphism of scheme sites.
If $f \sim g$, then $\widehat{f} \sim \widehat{g}$.
\end{lem}
\begin{proof}
Condition~\eqref{ssm:comp} is clear: if $(a,y) \leq (a_1,y_1)$, then $f^\sharp_{b,a_1}(y_1) = f^\sharp_{b,a}(\mathcal{O}^Y_{a,a_1}(y_1))$ divides $f^\sharp_{b,a}(y)$, so it is invertible whenever $f^\sharp_{b,a}(y)$ is.
Condition~\eqref{ssm:nat} follows from the uniqueness part of the universal property of localizations.
Condition~\eqref{ssm:top} follows from the corresponding condition for $f$ since $b \leq_f a$ implies $b \leq_{\widehat{f}} (a,1)$.

Let us verify condition~\eqref{ssm:meet}.
Suppose that $b \leq_{\widehat{f}} (a_1,y_1)$ and $b \leq_{\widehat{f}} (a_2,y_2)$.
By condition~\eqref{ssm:meet} for $f$ and \rlem{cover-prop}.\eqref{cover:restrict}, $b$ is covered by elements $b' \leq b$ such that $b' \leq_f a$ for some $a \leq a_1, a_2$.
For such $b'$, we have $b' \leq_{\widehat{f}} (a, \mathcal{O}^Y_{a,a_1}(y_1) \mathcal{O}^Y_{a,a_2}(y_2))$ and $(a, \mathcal{O}^Y_{a,a_1}(y_1) \mathcal{O}^Y_{a,a_2}(y_2)) \leq (a_i,y_i)$ for $i = 1,2$.

Let us verify condition~\eqref{ssm:loc}.
Let $(a',y') \leq (a,y)$ and $b \leq_{\widehat{f}} (a,y)$.
By \axref{SR}, $\mathcal{O}^Y_{a',a}$ is a localization at some $z \in \mathcal{O}^Y_a$.
By \rdefn{ring-loc}.\eqref{loc:surj}, there exist $w \in \mathcal{O}^Y_a$ and $c \in S_z$ such that $y' \cdot \mathcal{O}^Y_{a',a}(c) = \mathcal{O}^Y_{a',a}(w)$.
By \rlem{sat-loc}, $\mathcal{O}^{\mathrm{Sat}(Y)}_{(a',y'),(a,y)}$ is a localization at $z w / 1$.
By \rlem{loc-forall}, $b$ is covered by the set of elements $e \leq b$ for which there exists $e' \leq e$ such that $e' \leq_f a'$ and $\mathcal{O}^X_{e',e}$ is a localization at $f^\sharp_{e,a}(z)$.
We claim that every such $e$ belongs to $L_{\widehat{f}}((a,y),(a',y'))$.
Since $e \leq b$, we have $e \leq_{\widehat{f}} (a,y)$.
By \axref{SE}, there exists $e'' \leq e'$ such that $\mathcal{O}^X_{e'',e'}$ is a localization at $f^\sharp_{e',a'}(\mathcal{O}^Y_{a',a}(w)) = \mathcal{O}^X_{e',e}(f^\sharp_{e,a}(w))$.
By \rlem{loc-comp}, the map $\mathcal{O}^X_{e'',e}$ is a localization at $f^\sharp_{e,a}(z w) = \widehat{f}^\sharp_{e,(a,y)}(z w / 1)$.
Finally, $e'' \leq_{\widehat{f}} (a',y')$ since $e'' \leq_f a'$ and $f^\sharp_{e'',a'}(y')$ divides $f^\sharp_{e'',a'}(\mathcal{O}^Y_{a',a}(w))$, which is invertible.

Now, suppose that $f \sim g$, and let $b \leq_{\widehat{f}} (a,y)$.
Then $b \cover \{ b' \mid b' \leq_g a \}$, so, by \rlem{cover-prop}.\eqref{cover:restrict} and condition~\eqref{ssm:comp}, $b$ is covered by elements $b' \leq b$ such that $b' \leq_g a$.
For every such $b'$, we have $b' \leq_f a$, so $g^\sharp_{b',a}(y) = f^\sharp_{b',a}(y) = \mathcal{O}^X_{b',b}(f^\sharp_{b,a}(y))$ is invertible.
Thus, $b' \leq_{\widehat{g}} (a,y)$.
Finally, if $b \leq_{\widehat{f}} (a,y)$ and $b \leq_{\widehat{g}} (a,y)$, then $\widehat{f}^\sharp_{b,(a,y)} = \widehat{g}^\sharp_{b,(a,y)}$ since both extend $f^\sharp_{b,a} = g^\sharp_{b,a}$.
\end{proof}

\begin{lem}[counit-lift]
Let $X$ be a saturated scheme site, and let $f : X \to Y$ be a premorphism.
Then $\varepsilon_Y \circ \widehat{f} = f$.
\end{lem}
\begin{proof}
If $b \leq_{\widehat{f}} (a',y') \leq_{\varepsilon_Y} a$, then $b \leq_f a' \leq a$, so $b \leq_f a$.
Conversely, if $b \leq_f a$, then $b \leq_{\widehat{f}} (a,1) \leq_{\varepsilon_Y} a$.
Finally, the homomorphisms coincide since $\widehat{f}^\sharp_{b,(a',y')} \circ l_{y'} \circ \mathcal{O}^Y_{a',a} = f^\sharp_{b,a'} \circ \mathcal{O}^Y_{a',a} = f^\sharp_{b,a}$.
\end{proof}

\begin{lem}[counit-mono]
Let $X$ and $Y$ be scheme sites, and let $f,g : X \to \mathrm{Sat}(Y)$ be premorphisms such that $\varepsilon_Y \circ f \sim \varepsilon_Y \circ g$.
Then $f \sim g$.
\end{lem}
\begin{proof}
Note that $(a',y') \leq (a,1)$ whenever $a' \leq a$, and that $l_1$ is an isomorphism.
Thus, if $b \leq_f (a,1)$ and $b \leq_g (a,1)$, then $b \leq_{\varepsilon_Y \circ f} a$, $b \leq_{\varepsilon_Y \circ g} a$, and
\[ f^\sharp_{b,(a,1)} \circ l_1 = (\varepsilon_Y \circ f)^\sharp_{b,a} = (\varepsilon_Y \circ g)^\sharp_{b,a} = g^\sharp_{b,(a,1)} \circ l_1, \]
so $f^\sharp_{b,(a,1)} = g^\sharp_{b,(a,1)}$.

Let $b \leq_f (a,y)$.
Then $b \leq_{\varepsilon_Y \circ f} a$, so $b$ is covered by the set of elements $b'$ such that $b' \leq_{\varepsilon_Y \circ g} a$.
Every such $b'$ satisfies $b' \leq_g (a,1)$.
By \rlem{cover-prop}.\eqref{cover:restrict} and condition~\eqref{ssm:comp}, $b$ is covered by the set of elements $b' \leq b$ such that $b' \leq_g (a,1)$.
By \rlem{cover-prop}.\eqref{cover:trans}, it is enough to show that every such $b'$ is covered by $\{ b'' \mid b'' \leq_g (a,y) \}$.
We have $b' \leq_f (a,y) \leq (a,1)$.
By \rlem{sat-loc}, $\mathcal{O}^{\mathrm{Sat}(Y)}_{(a,y),(a,1)}$ is a localization at $y/1$.
The element $g^\sharp_{b',(a,1)}(y/1) = f^\sharp_{b',(a,1)}(y/1) = f^\sharp_{b',(a,y)}(y/1)$ is invertible, so the required property follows from \rlem{restrict}.
The symmetric implication is proved in the same way.

Finally, if $b \leq_f (a,y)$ and $b \leq_g (a,y)$, then
\[ f^\sharp_{b,(a,y)} \circ \mathcal{O}^{\mathrm{Sat}(Y)}_{(a,y),(a,1)} = f^\sharp_{b,(a,1)} = g^\sharp_{b,(a,1)} = g^\sharp_{b,(a,y)} \circ \mathcal{O}^{\mathrm{Sat}(Y)}_{(a,y),(a,1)}. \]
Since localizations are epimorphisms, we have $f^\sharp_{b,(a,y)} = g^\sharp_{b,(a,y)}$.
\end{proof}

\begin{thm}[coreflection]
For every saturated scheme site $X$ and every scheme site $Y$, the maps $f \mapsto \widehat{f}$ and $h \mapsto \varepsilon_Y \circ h$
are mutually inverse bijections between morphisms $X \to Y$ and morphisms $X \to \mathrm{Sat}(Y)$.
Consequently, $\cat{SatSchSite}$ is a coreflective subcategory of $\cat{SchSite}$:
the saturation extends to a functor $\mathrm{Sat} : \cat{SchSite} \to \cat{SatSchSite}$ right adjoint to the inclusion, with counit $\varepsilon$.
\end{thm}
\begin{proof}
By \rlem{lift}, the map $f \mapsto \widehat{f}$ is well-defined on equivalence classes.
By \rlem{counit-lift}, $\varepsilon_Y \circ \widehat{f} = f$.
Conversely, for every premorphism $h : X \to \mathrm{Sat}(Y)$, we have $\varepsilon_Y \circ \widehat{\varepsilon_Y \circ h} = \varepsilon_Y \circ h$ by \rlem{counit-lift},
so $\widehat{\varepsilon_Y \circ h} \sim h$ by \rlem{counit-mono}.
The bijection $h \mapsto \varepsilon_Y \circ h$ is clearly natural in $X$.
Thus, $\varepsilon_Y$ is a coreflection of $Y$ into $\cat{SatSchSite}$.
\end{proof}

Explicitly, the action of the functor $\mathrm{Sat}$ on a morphism $f : X \to Y$ is given by $\mathrm{Sat}(f) = \widehat{f \circ \varepsilon_X}$.

\section{Equivalence with schemes}
\label{sec:equivalence}

In this section, we prove that the category $\cat{SatSchSite}$ of saturated scheme sites is equivalent to the category of schemes.
Together with \rthm{coreflection}, this implies that the category of schemes is equivalent to a coreflective subcategory of $\cat{SchSite}$.
Since we work constructively, we define schemes as locally ringed locales instead of locally ringed spaces.
Classically, these approaches are equivalent, but constructively only the localic approach works well.

We recall necessary definitions first.
For background on locales and sheaves on them, see \cite[Chapter~II]{stone-spaces}, \cite[Part~C]{elephant}, and \cite{maclane-moerdijk}.
A \emph{frame} is a complete distributive lattice.
A map of frames is a function that preserves finite meets and arbitrary joins.
The category of locales is the opposite of the category of frames.
Given a locale $L$, we denote its underlying frame by $O_L$.
Elements of $O_L$ are called \emph{opens} of $L$.
For a map of locales $f : L \to M$, the corresponding frame map is written as $f^* : O_M \to O_L$.

A \emph{sheaf} of rings on a locale $L$ is a presheaf $\mathcal{O}^L$ on $O_L$ that satisfies the usual sheaf condition.
A \emph{ringed locale} is a pair $(L,\mathcal{O}^L)$ consisting of a locale and a sheaf of rings on it.
A \emph{locally ringed locale} is a ringed locale $(L,\mathcal{O}^L)$ that satisfies the following condition:
for every open $a \in O_L$ and every sequence of elements $x_1, \ldots, x_n \in \mathcal{O}^L_a$ such that $x_1 + \ldots + x_n = 1$, we have
\[ a \leq \bigvee \{ b \in O_L \mid b \leq a, \exists i, \mathcal{O}^L_{b,a}(x_i) \text{ is invertible} \}. \]

\begin{defn}[lrl-map]
A \emph{morphism of ringed locales} $(L,\mathcal{O}^L)$ and $(M,\mathcal{O}^M)$ is a pair $(f,f^\sharp)$ consisting of a map $f : L \to M$ and a natural transformation $f^\sharp_a : \mathcal{O}^M_a \to \mathcal{O}^L_{f^*(a)}$.
A \emph{morphism of locally ringed locales} $(L,\mathcal{O}^L)$ and $(M,\mathcal{O}^M)$ is a morphism of ringed locales $(f,f^\sharp)$ that satisfies the following condition:
for all $b \in O_L$, $a \in O_M$, and $x \in \mathcal{O}^M_a$, such that $b \leq f^*(a)$ and $\mathcal{O}^L_{b,f^*(a)}(f^\sharp_a(x))$ is invertible, we have
\[ b \leq \bigvee \{ f^*(a') \mid a' \leq a, \mathcal{O}^M_{a',a}(x) \text{ is invertible} \}. \]
We denote the category of locally ringed locales by $\cat{LRL}$.
\end{defn}

Locally ringed locales (and, more generally, locally ringed toposes) go back to \cite{hakim}; see also \cite{blechschmidt}.

Schemes are usually defined as locally ringed locales which are locally affine.
The condition of being locally affine is typically formulated in terms of spectra of rings.
However, for our purposes it will be more convenient to use a more explicit formulation.
We will show later that this formulation is equivalent to the usual one.

Let $(L,\mathcal{O}^L)$ be a ringed locale.
Let $a \in \mathcal{O}_L$ be an open, and let $x \in \mathcal{O}^L_a$ be an element.
Define $P(a,x) \in O_L$ as follows:
\[ P(a,x) = \bigvee \{ b \in O_L \mid b \leq a, \mathcal{O}^L_{b,a}(x) \text{ is invertible} \}. \]
Then we have $P(a,x) \leq a$.
Also, since $\mathcal{O}^L$ is a sheaf, $\mathcal{O}^L_{P(a,x),a}(x)$ is invertible.
Thus, $P(a,x)$ is the largest open $b$ such that $b \leq a$ and $\mathcal{O}^L_{b,a}(x)$ is invertible.
We call opens of the form $P(a,x)$ \emph{principal opens}.

\begin{defn}[affine-open]
We say that an open $a \in \mathcal{O}_L$ is an \emph{affine open} if the following conditions are satisfied:
\begin{enumerate}
\item \label{affine:loc} For every $x \in \mathcal{O}^L_a$, the map $\mathcal{O}^L_{P(a,x),a}$ is a localization at $x$.
\item \label{affine:surj} For every $b \leq a$, we have
\[ b \leq \bigvee \{ c \mid c \leq b, \exists x \in \mathcal{O}^L_a, c = P(a,x) \}. \]
\item \label{affine:cover} For all $x \in \mathcal{O}^L_a$ and $U \subseteq O_L$ such that $P(a,x) \leq \bigvee U$,
the ideal of $\mathcal{O}^L_{P(a,x)}$ generated by the following set is the unit ideal:
\[ \{ y \mid \exists c \in U, \exists b, b \leq P(a,x), b \leq c, \mathcal{O}^L_{b,P(a,x)} \text{ is a localization at } y \}. \]
\end{enumerate}
\end{defn}

A \emph{scheme} is a locally ringed locale which can be covered by affine opens.
This means that $\bigvee \{ a \in O_L \mid a \text{ is an affine open} \}$ is the greatest element of $O_L$.
Morphisms of schemes are just morphisms of underlying locally ringed locales.
We denote the category of schemes by $\cat{Sch}$.

\subsection{From scheme sites to schemes}
\label{sec:scheme-site-to-scheme}

We now turn to the construction of a scheme from a saturated scheme site.
As a first step, we need to show that the structure presheaf is a sheaf with respect to the covering relation $\cover$.
For this, we will use the following standard algebraic lemma.

\begin{lem}[loc-amalg]
Let $R$ be a ring, and let $\{b_i\}_{i \in I}$ be a family of elements of $R$ such that the ideal generated by the $b_i$ is the unit ideal.
For each $i \in I$, let $x_i \in R[b_i^{-1}]$ be an element, and assume that for all $i,j \in I$, the images of $x_i$ and $x_j$ in $R[(b_i b_j)^{-1}]$ coincide.
Then there exists a unique element $x \in R$ such that for every $i \in I$, the image of $x$ in $R[b_i^{-1}]$ is equal to $x_i$.
\end{lem}

For simplicity, we have formulated this lemma for standard localizations of the form $R[b_i^{-1}]$.
However, it remains valid for arbitrary localizations in the sense of \rdefn{ring-loc}.

\begin{lem}[str-sheaf]
Let $X$ be a scheme site.
Then the presheaf $\mathcal{O}^X$ is a sheaf with respect to the covering relation $\cover$.
More precisely, for every $a \in X$ and every covering $a \cover U$ such that elements of $U$ are less than or equal to $a$, the following holds.

If for each $b \in U$ we are given an element $x_b \in \mathcal{O}^X_b$ such that
for all $b_1,b_2 \in U$ and all $c \leq b_1,b_2$ we have $\mathcal{O}^X_{c,b_1}(x_{b_1}) = \mathcal{O}^X_{c,b_2}(x_{b_2})$,
then there exists a unique element $x \in \mathcal{O}^X_a$ such that for every $b \in U$ we have $\mathcal{O}^X_{b,a}(x) = x_b$.
\end{lem}
\begin{proof}
To apply \rlem{loc-amalg}, we need to show that the family $\{ x_b \}_{b \in U}$ satisfies the required compatibility condition.
To do this, it is enough to check that localizations $\mathcal{O}^X_a[(b_1 b_2)^{-1}]$ can be represented as $\mathcal{O}^X_c$ for some $c \leq b_1,b_2$.
This follows directly from \axref{SS}, which ensures that such a $c$ exists and that $\mathcal{O}^X_c$ is obtained from $\mathcal{O}^X_a$ by localizing at $b_1 b_2$.
\end{proof}

Now, given a scheme site $X$, we construct the corresponding locale and a structure sheaf on it.
In the terminology of \cite{vickers-compact}, the preorder $X$ together with the covering relation $\cover$ forms a flat site $(X,\cover)$.
Any such site generates a locale, which we denote by $\widetilde{X}$.
The elements of $\widetilde{X}$ are $\cover$-ideals, that is, subsets $U \subseteq X$ such that
\[
a \cover U \Rightarrow a \in U.
\]

For every element $a \in X$, we write $\overline{a}$ for the $\cover$-ideal $\{ x \in X \mid x \cover \{a\} \}$.
We call elements of the form $\overline{a}$ \emph{basic opens}.
The map $a \mapsto \overline{a} : X \to \widetilde{X}$ is a dense morphism of sites (see, for example, \cite[Section~11]{shulman-flat} or \cite[Section~9.4]{minichiello} for the definition and properties of dense morphisms of sites).
In particular, it induces an equivalence between the category of sheaves on $\widetilde{X}$ and the category of sheaves on $X$ (this is an instance of the comparison lemma; see \cite[Theorem~C2.2.3]{elephant}).
Consequently, every sheaf $\mathcal{F}$ on $X$ uniquely extends to a sheaf $\widetilde{\mathcal{F}}$ on $\widetilde{X}$:
\[ \xymatrix{ X \ar[r]^-{\mathcal{F}} \ar[d] & \cat{CRing} \\
              \widetilde{X} \ar@{-->}[ur]_{\widetilde{\mathcal{F}}}
            } \]
This extension can be explicitly defined as follows:
\[ \widetilde{\mathcal{F}}_U = \lim_{b \in U} \mathcal{F}_b. \]
Given an element $b \in U$, we write $\pi_b$ for the canonical projection $\widetilde{\mathcal{F}}_U \to \mathcal{F}_b$.

Thus, starting from a scheme site $(X,\mathcal{O}^X)$, we obtain a ringed locale $(\widetilde{X},\widetilde{\mathcal{O}^X})$.
We now show that it is a scheme.

\begin{lem}
For every saturated scheme site $(X,\mathcal{O}^X)$, the ringed locale $(\widetilde{X},\widetilde{\mathcal{O}^X})$ is locally ringed.
\end{lem}
\begin{proof}
We verify this property only for basic opens since the generalization to arbitrary opens is straightforward.
We need to show that for every $a \in X$ and every sequence of elements $x_1, \ldots, x_n \in \mathcal{O}^X_a$ such that $x_1 + \ldots + x_n = 1$, the element $a$ is covered by the following set:
\[ U = \{ b \in X \mid b \leq a, \exists i, \mathcal{O}^X_{b,a}(x_i) \text{ is invertible} \}. \]

It is enough to show that for every $i$, the element $x_i$ belongs to $C(a,U)$.
This means we need to show that there exists $b \in X$ such that $\mathcal{O}^X_{b,a}$ is a localization at $x_i$, but this is exactly \axref{SE}.
\end{proof}

\begin{lem}[po-char]
Let $a,b \in X$ and $x \in \widetilde{\mathcal{O}^X}_{\overline{a}}$ be such that $b \leq a$ and $\mathcal{O}^X_{b,a}$ is a localization at $\pi_a(x)$.
Then $P(\overline{a},x) = \overline{b}$.
\end{lem}
\begin{proof}
The map $\pi_b : \widetilde{\mathcal{O}^X}_{\overline{b}} \to \mathcal{O}^X_b$ is an isomorphism by construction of $\widetilde{\mathcal{O}^X}$.
It follows that $\widetilde{\mathcal{O}^X}_{\overline{b},\overline{a}}(x)$ is invertible.
Hence $\overline{b} \leq P(\overline{a},x)$.

To prove that $P(\overline{a},x) \leq \overline{b}$, we need to show that for every $c \in X$
such that $c \leq a$ and $\mathcal{O}^X_{c,a}(\pi_a(x))$ is invertible, we have $c \cover \{ b \}$.
This follows from \rlem{inv-cover}.
\end{proof}

\begin{prop}
For every saturated scheme site $(X,\mathcal{O}^X)$, the locally ringed locale $(\widetilde{X},\widetilde{\mathcal{O}^X})$ is a scheme.
\end{prop}
\begin{proof}
It is enough to show that basic opens are affine since every open is a join of basic opens.
Let $a \in X$ be an element.
We verify the conditions of \rdefn{affine-open} for $\overline{a}$.

\eqref{affine:loc}
We need to show that for every $x \in \widetilde{\mathcal{O}^X}_{\overline{a}}$,
the map $\widetilde{\mathcal{O}^X}_{P(\overline{a},x),\overline{a}}$ is a localization at $x$.
Define $y \in \mathcal{O}^X_a$ as $\pi_a(x)$.
By \axref{SE}, there exists $b \in X$ such that $b \leq a$ and $\mathcal{O}^X_{b,a}$ is a localization at $y$.
Consider the following diagram:
\[ \xymatrix{ \widetilde{\mathcal{O}^X}_{\overline{a}} \ar[r]^{\pi_a} \ar[d]_{\widetilde{\mathcal{O}^X}_{P(\overline{a},x),\overline{a}}} & \mathcal{O}^X_a \ar[d]^{\mathcal{O}^X_{b,a}} \\
              \widetilde{\mathcal{O}^X}_{P(\overline{a},x)} \ar[r]_-{\pi_b}                                                               & \mathcal{O}^X_b
            } \]

The map $\pi_a : \widetilde{\mathcal{O}^X}_{\overline{a}} \to \mathcal{O}^X_a$ is an isomorphism by construction of $\widetilde{\mathcal{O}^X}$.
By \rlem{po-char}, the map $\pi_b : \widetilde{\mathcal{O}^X}_{P(\overline{a},x)} \to \mathcal{O}^X_b$ is also an isomorphism.
It follows that $\widetilde{\mathcal{O}^X}_{P(\overline{a},x),\overline{a}}$ is a localization at $x$.

\eqref{affine:surj}
It is enough to prove this property for basic opens $\overline{b} \leq \overline{a}$ such that $b \leq a$,
but for such an open, it follows from \axref{SR} and \rlem{po-char}.

\eqref{affine:cover}
Let $x \in \widetilde{\mathcal{O}^X}_{\overline{a}}$ and $U \subseteq \widetilde{X}$ be such that $P(\overline{a},x) \leq \bigvee U$.
Let $b \leq a$ be such that $\mathcal{O}^X_{b,a}$ is a localization at $\pi_a(x)$.
By \rlem{po-char}, we have $P(\overline{a},x) = \overline{b}$.
It follows that $b \cover \bigcup U$.

Applying the isomorphism $\pi_b : \widetilde{\mathcal{O}^X}_{P(\overline{a},x)} \to \mathcal{O}^X_b$
to the set described in the condition produces exactly the set $C(b, \bigcup U)$, which generates the unit ideal.
This concludes the proof.
\end{proof}

\subsection{Affine opens}
\label{sec:affine-opens}

Affine opens are usually defined in terms of spectra of rings.
We now show that our definition is equivalent to the standard one.
To do this, we first recall a few standard definitions and lemmas.

The \emph{spectrum} of a ring $R$ is the scheme $\widetilde{\mathrm{Spec}(R)}$ associated to the saturated scheme site $\mathrm{Spec}(R) = \mathrm{Sat}(\mathrm{pt}_R)$ defined in \rexample{spec}.
This construction gives us a functor $\mathrm{Spec} : \cat{CRing}^{\mathrm{op}} \to \cat{Sch}$.
This functor has a left adjoint $\Gamma : \cat{Sch} \to \cat{CRing}^{\mathrm{op}}$ defined by $\Gamma(X) = \mathcal{O}^X_\top$,
where $\top$ is the largest open of the underlying locale.

We will need an explicit description of the unit of the adjunction $\eta_X : X \to \mathrm{Spec}(\Gamma(X))$.
For an open $U$ of $\mathrm{Spec}(\Gamma(X))$, we define $\eta_X^*(U)$ as follows:
\[ \eta_X^*(U) = \bigvee_{x \in U} P(\top,x). \]
We now define $\eta_X^\sharp : \mathcal{O}^{\mathrm{Spec}(\Gamma(X))}_U \to \mathcal{O}^X_{\eta_X^*(U)}$.
We have the following isomorphisms:
\begin{align*}
\mathcal{O}^{\mathrm{Spec}(\Gamma(X))}_U & \simeq \lim_{x \in U} \mathcal{O}^X_\top[x^{-1}], \\
\mathcal{O}^X_{\eta_X^*(U)} & \simeq \lim_{x \in U} \mathcal{O}^X_{P(\top,x)}.
\end{align*}
The first isomorphism holds by the definition of $\mathrm{Spec}$.
The second one follows from the sheaf condition.
Thus, we can define $\eta_X^\sharp$ by the universal property of limits as the unique map such that the following diagram commutes for every $x \in U$:
\[ \xymatrix{ \mathcal{O}^{\mathrm{Spec}(\Gamma(X))}_U \ar[r]^-{\eta_X^\sharp} \ar[d]_{\pi_x} & \mathcal{O}^X_{\eta_X^*(U)} \ar[d]^{\mathcal{O}^X_{P(\top,x), \eta_X^*(U)}} \\
              \mathcal{O}^X_\top[x^{-1}] \ar[r] & \mathcal{O}^X_{P(\top,x)}
            } \]
The bottom map is defined by the universal property of $\mathcal{O}^X_\top[x^{-1}]$.

The following lemma shows that $P(\top,-)$ preserves coverings:

\begin{lem}[eta-cover]
Let $X$ be a locally ringed locale.
For all $y \in \Gamma(X)$ and $U \subseteq \Gamma(X)$ such that $y \cover U$, we have
\[ P(\top,y) \leq \bigvee_{x \in U} P(\top,x). \]
\end{lem}
\begin{proof}
By assumption, there exist $x_1, \ldots, x_n \in U$ such that a power of $y$ is a linear combination of $x_1, \ldots, x_n$.
Let $b \in O_X$ be an open such that $\mathcal{O}^X_{b,\top}(y)$ is invertible.
We need to show that $b \leq \bigvee_{x \in U} P(\top,x)$.
Since $\mathcal{O}^X_{b,\top}(y)$ is invertible, there exist $\alpha_1, \ldots, \alpha_n \in \mathcal{O}^X_b$ such that
\[ \alpha_1 \mathcal{O}^X_{b,\top}(x_1) + \ldots + \alpha_n \mathcal{O}^X_{b,\top}(x_n) = 1. \]
Since $X$ is locally ringed, we have
\[ b \leq \bigvee \{ c \mid c \leq b, \exists i, \mathcal{O}^X_{c,b}(\alpha_i \mathcal{O}^X_{b,\top}(x_i)) \text{ is invertible} \} \leq \bigvee_{x \in U} P(\top,x). \]
\end{proof}

We are now ready to prove that our definition of affine opens is equivalent to the standard one:

\begin{prop}
An open of a scheme $a \in O_X$ is affine if and only if the restriction of $X$ to $a$ is isomorphic to a spectrum of a ring.
\end{prop}
\begin{proof}
It is easy to see that $a$ is an affine open if and only if $\top$ is an affine open in the $X$ restricted to $a$.
Thus, we may assume that $a = \top$.

Clearly, the spectrum of a ring satisfies the conditions of \rdefn{affine-open}.
We prove the converse.
Assume that $\top$ is an affine open.
The opposite of the category of rings is a reflective subcategory of $\cat{LRL}$.
It follows that a locally ringed locale $X$ is isomorphic to a spectrum of a ring if and only if $\eta_X$ is an isomorphism.
To prove that $\eta_X$ is an isomorphism, it is enough to show that $\eta_X^*$ and $\eta_X^\sharp$ are both bijective.

By \eqref{affine:loc}, the map $\mathcal{O}^X_\top[x^{-1}] \to \mathcal{O}^X_{P(\top,x)}$ in the definition of $\eta_X^\sharp$ is an isomorphism.
It follows that $\eta_X^\sharp$ is also an isomorphism.

The surjectivity of $\eta_X^*$ follows from \eqref{affine:surj}.
Indeed, let $b$ be an open of $X$.
Define $U$ as $\{ y \mid y \cover \{ x \mid P(\top,x) \leq b \}$ \}.
Then \rlem{eta-cover} implies that $\eta_X^*(U) \leq b$ and $\eqref{affine:surj}$ implies that $b \leq \eta_X^*(U)$.

Finally, to prove injectivity of $\eta_X^*$, it is enough to show that $\eta_X^*(V) \leq \eta_X^*(U)$ implies that $V \subseteq U$.
Let $x$ be an element in $V$.
We need to show that $x \cover U$.

By assumption,
\[ P(\top,x) \leq \bigvee_{y \in U} P(\top,y). \]
By \eqref{affine:cover}, the following set generates the unit ideal of $\mathcal{O}^X_{P(\top,x)}$:
\[ \{ z \mid \exists y \in U, \exists b, b \leq P(\top,x), b \leq P(\top,y), \mathcal{O}^X_{b,P(\top,x)} \text{ is a localization at } z \}. \]
It follows that a power of $x$ belongs to the ideal generated by the following set:
\[ \{ z \in \mathcal{O}^X_\top \mid \exists y \in U, \exists b, b \leq P(\top,x), b \leq P(\top,y), \mathcal{O}^X_{b,\top} \text{ is a localization at } z \}. \]
Thus, it is enough to show that every element of this set is covered by $U$.
If we have $y \in U$ such that $b \leq P(\top,y)$, then $\mathcal{O}^X_{b,\top}(y)$ is invertible.
It follows that $y$ divides a power of $z$, which implies $z \cover \{ y \}$.
This completes the proof.
\end{proof}

\subsection{Functoriality}
\label{sec:functor}

In Section~\ref{sec:scheme-site-to-scheme}, we constructed a scheme $\widetilde{X}$ from a saturated scheme site $X$.
Now, we show that this construction is functorial.

To define this functor, it is convenient to introduce an auxiliary notion of a presentation morphism from a locally ringed locale to a scheme site:

\begin{defn}[pres-map]
Let $(L,\mathcal{O}^L)$ be a locally ringed locale, and let $Y$ be a saturated scheme site.
A \emph{presentation morphism} from $(L,\mathcal{O}^L)$ to $Y$ is a pair $(f^*,f^\sharp)$ consisting of a monotone map $f^* : Y \to O_L$
and a natural transformation $f^\sharp_a : \mathcal{O}^Y_a \to \mathcal{O}^L_{f^*(a)}$ satisfying the following conditions:
\begin{enumerate}
\item \label{pres:top} Locale $L$ is covered by the image of $f^*$, that is, $\top \leq \bigvee_{a \in Y} f^*(a)$.
\item \label{pres:meet} For all $a,b \in Y$, we have $f^*(a) \wedge f^*(b) \leq \bigvee \{ f^*(c) \mid c \leq a, c \leq b \}$.
\item \label{pres:loc} For all $a,a' \in Y$, $x \in \mathcal{O}^Y_a$, and $b \in O_L$, if $a' \leq a$, $\mathcal{O}^Y_{a',a}$ is a localization at $x$,
$b \leq f^*(a)$, and $\mathcal{O}^L_{b,f^*(a)}(f^\sharp_a(x))$ is invertible, then $b \leq f^*(a')$.
\end{enumerate}
\end{defn}

We now show that there is a bijection between presentation morphisms from $L$ to $Y$ and morphisms of locally ringed locales between $L$ and $\widetilde{Y}$.
Recall that, for every $a \in Y$, the projection map $\pi_a : \widetilde{\mathcal{O}^Y}_{\overline{a}} \to \mathcal{O}^Y_a$ is an isomorphism.
Now, given a morphism of locally ringed locales $f : L \to \widetilde{Y}$, we define the corresponding presentation morphism $\psi'(f) : L \to Y$ as follows:
\begin{align*}
\psi'(f)^*(a) & = f^*(\overline{a}), \\
\psi'(f)^\sharp_a & = f^\sharp_{\overline{a}} \circ \pi_a^{-1}.
\end{align*}
The first two conditions of \rdefn{pres-map} follow from the fact that $f^*$ preserves finite meets.
The last one follows from the fact that $f$ satisfies the condition of \rdefn{lrl-map} and \rlem{po-char}.

Conversely, given a presentation morphism $f : L \to Y$, we define the corresponding morphism of locally ringed locales $\psi(f) : L \to \widetilde{Y}$.
For every $U \in O_{\widetilde{Y}}$, define $\psi(f)^*(U)$ as $\bigvee_{a \in U} f^*(a)$.
By the sheaf condition, we have
\begin{align*}
\widetilde{\mathcal{O}^Y}_U & \simeq \lim_{a \in U} \mathcal{O}^Y_a, \\
\mathcal{O}^L_{\psi(f)^*(U)} & \simeq \lim_{a \in U} \mathcal{O}^L_{f^*(a)}.
\end{align*}
Thus, we can define $\psi(f)^\sharp$ as the unique arrow that fits in the following diagram for all $a \in U$:
\[ \xymatrix{ \widetilde{\mathcal{O}^Y}_U \ar[r]^-{\psi(f)^\sharp_U} \ar[d]_{\pi_a} & \mathcal{O}^L_{\psi(f)^*(U)} \ar[d]^{\mathcal{O}^L_{f^*(a),\psi(f)^*(U)}} \\
              \mathcal{O}^Y_a \ar[r]_-{f^\sharp_a}                                  & \mathcal{O}^L_{f^*(a)}
            } \]

\begin{lem}
For every presentation morphism $f : L \to Y$, the construction of $\psi(f)$ is correct.
\end{lem}
\begin{proof}
First, we verify that $\psi(f)^*$ is a morphism of locales.
Conditions~\eqref{pres:top} and \eqref{pres:meet} of \rdefn{pres-map} guarantee that $f^*$ is a flat functor.
Thus, to prove that it induces a morphism of locales, we just need to verify that it respects covers.

Let $a \in Y$ and $U \subseteq Y$ be such that $a \cover U$.
We need to prove that $f^*(a) \leq \bigvee_{c \in U} f^*(c)$.
By assumption, there exists a sequence $x_1, \ldots, x_n \in C(a,U)$ such that $1$ is a linear combination of $x_1, \ldots, x_n$.
It follows that $1$ is a linear combination of $f^\sharp_a(x_1), \ldots, f^\sharp_a(x_n)$ in $\mathcal{O}^L_{f^*(a)}$.
Since $L$ is locally ringed, we have
\[ f^*(a) \leq \bigvee \{ b \mid b \leq f^*(a), \exists i, \mathcal{O}^L_{b,f^*(a)}(f^\sharp_a(x_i)) \text{ is invertible} \}. \]

Let $b \in O_L$ be such that $b \leq f^*(a)$ and $\mathcal{O}^L_{b,f^*(a)}(f^\sharp_a(x_i))$ is invertible.
We need to show that $b \leq \bigvee_{c \in U} f^*(c)$.
Since $x_i \in C(a,U)$, there exist $a' \in Y$ and $c \in U$ such that $a' \leq a$, $a' \leq c$, and $\mathcal{O}^Y_{a',a}$ is a localization at $x_i$.
By \rdefn{pres-map}.\eqref{pres:loc}, we know that $b \leq f^*(a')$.
It follows that $b \leq f^*(a') \leq f^*(c) \leq \bigvee_{c \in U} f^*(c)$.

The fact that $\psi(f)$ is a morphism of locally ringed locales follows from \rdefn{pres-map}.\eqref{pres:loc} and \axref{SE}.
\end{proof}

\begin{lem}[psi-bij]
The above constructions define a bijection between morphisms of locally ringed locales $L \to \widetilde{Y}$ and presentation morphisms $L \to Y$.
\end{lem}
\begin{proof}
It is clear that $\psi'(\psi(f)) = f$ for every $f : L \to Y$.
The other direction follows from the facts that $Y$ is a dense subsite of $\widetilde{Y}$ and that the extension of a sheaf from a dense subsite to the whole site is unique.
\end{proof}

We now construct a bijection between morphisms of saturated scheme sites $X \to Y$ and presentation morphisms $\widetilde{X} \to Y$.
To do this, it is convenient to introduce an alternative description of morphisms of scheme sites.

\begin{defn}[complete]
A premorphism $f : X \to Y$ is \emph{complete} if, for all $b \in X$ and $a \in Y$,
\[ b \cover \{ b' \mid b' \leq_f a \} \Rightarrow b \leq_f a. \]
\end{defn}

\begin{lem}[closure]
Let $f : X \to Y$ be a premorphism.
Then there is a complete premorphism $\overline{f}$ such that
\[ b \leq_{\overline{f}} a \text{ if and only if } b \cover \{ b' \mid b' \leq_f a \}, \]
and $\overline{f}^\sharp_{b,a}$ is the unique map fitting in the following diagram for every $b'$ such that $b' \leq b$ and $b' \leq_f a$:
\[ \xymatrix{ \mathcal{O}^Y_a \ar[r]^{\overline{f}^\sharp_{b,a}} \ar[dr]_{f^\sharp_{b',a}} & \mathcal{O}^X_b \ar[d]^{\mathcal{O}^X_{b',b}} \\
                                                                                           & \mathcal{O}^X_{b'}
            } \]
Moreover, $f \sim \overline{f}$, and equivalent complete premorphisms are equal.
\end{lem}
\begin{proof}
First, let us show that $\overline{f}^\sharp_{b,a}$ is well-defined.
By \rlem{cover-prop}.\eqref{cover:restrict} and condition~\eqref{ssm:comp}, $b$ is covered by the set of elements $b' \leq b$ such that $b' \leq_f a$.
For every $x \in \mathcal{O}^Y_a$, the family $f^\sharp_{b',a}(x)$ is compatible by condition~\eqref{ssm:nat}.
By \rlem{str-sheaf}, there exists a unique element $\overline{f}^\sharp_{b,a}(x)$ making the diagram commute.
Uniqueness implies that $\overline{f}^\sharp_{b,a}$ is a ring homomorphism.
If $b \leq_f a$, then $f^\sharp_{b,a}$ satisfies the same condition by \eqref{ssm:nat}, so $\overline{f}^\sharp_{b,a} = f^\sharp_{b,a}$.
In particular, $\leq_f$ is contained in $\leq_{\overline{f}}$, and $f \sim \overline{f}$.

Now, we verify the conditions of \rdefn{ssm}.
Condition~\eqref{ssm:comp} follows from \rlem{cover-prop}.\eqref{cover:left} and \eqref{cover:mono},
and completeness of $\overline{f}$ follows from \rlem{cover-prop}.\eqref{cover:trans}.
Condition~\eqref{ssm:nat} follows from \rlem{str-sep}.
Condition~\eqref{ssm:top} holds since $\leq_f$ is contained in $\leq_{\overline{f}}$.

Let us verify condition~\eqref{ssm:meet}.
Suppose that $b \leq_{\overline{f}} a_1$ and $b \leq_{\overline{f}} a_2$.
By \rlem{cover-prop}.\eqref{cover:meet} and condition~\eqref{ssm:comp} for $f$, we have $b \cover \{ c \mid c \leq_f a_1, c \leq_f a_2 \}$.
By condition~\eqref{ssm:meet} for $f$, every such $c$ is covered by $\{ c' \mid \exists a \in Y, a \leq a_1, a \leq a_2, c' \leq_f a \}$.
Since $\leq_f$ is contained in $\leq_{\overline{f}}$, the required property follows from \rlem{cover-prop}.\eqref{cover:trans}.

Let us verify condition~\eqref{ssm:loc}.
Suppose that $b \leq_{\overline{f}} a$ and $a' \leq a$.
Then $b \cover \{ b_1 \mid b_1 \leq_f a \}$, and every such $b_1$ is covered by $L_f(a,a')$.
Since $\leq_f$ is contained in $\leq_{\overline{f}}$ and the maps $\overline{f}^\sharp$ and $f^\sharp$ coincide on $\leq_f$,
we have $L_f(a,a') \subseteq L_{\overline{f}}(a,a')$.
Thus, the required property follows from \rlem{cover-prop}.\eqref{cover:trans}.

Finally, if $f$ and $g$ are complete and $f \sim g$, then $b \leq_f a$ implies $b \cover \{ b' \mid b' \leq_g a \}$, so $b \leq_g a$, and vice versa.
Thus, the relations $\leq_f$ and $\leq_g$ coincide, and $f^\sharp = g^\sharp$ by the definition of $\sim$.
\end{proof}

Thus, complete premorphisms are unique representatives of their equivalence classes.
Consequently, morphisms of scheme sites can be defined as complete premorphisms rather than equivalence classes of premorphisms.
In the rest of this section, we always represent morphisms by complete premorphisms.
Note that the complete representative of the composition of morphisms $f$ and $g$ is $\overline{g \circ f}$,
and the complete representative of the identity morphism is given by $b \leq_{\overline{\mathrm{id}_X}} a$ if and only if $b \cover \{ a \}$.

Given a complete premorphism of saturated scheme sites $f : X \to Y$, we define the corresponding presentation morphism $\varphi(f) : \widetilde{X} \to Y$.
For each $a \in Y$, we define $\varphi(f)^*(a)$ as $\{ b \mid b \leq_f a \}$.
This subset is a $\cover$-ideal since $f$ is complete.
The map $\varphi(f)^\sharp_a$ is defined as the unique map that fits in the following diagram for every $b$ such that $b \leq_f a$:
\[ \xymatrix{ \mathcal{O}^Y_a \ar[r]^-{\varphi(f)^\sharp_a} \ar[dr]_{f^\sharp_{b,a}} & \widetilde{\mathcal{O}^X}_{\varphi(f)^*(a)} \ar[d]^{\pi_b} \\
                                                                                     & \mathcal{O}^X_b
            } \]
Conditions~\eqref{pres:top} and \eqref{pres:meet} of \rdefn{pres-map} follow
from conditions~\eqref{ssm:top} and \eqref{ssm:meet} of \rdefn{ssm}, respectively.
To verify condition~\eqref{pres:loc}, let $a' \leq a$, let $\mathcal{O}^Y_{a',a}$ be a localization at $x$, and let $U \leq \varphi(f)^*(a)$ be an open such that the restriction of $\varphi(f)^\sharp_a(x)$ to $U$ is invertible.
Then for every $b \in U$, we have $b \leq_f a$ and $f^\sharp_{b,a}(x)$ is invertible.
By \rlem{restrict}, we have $b \cover \{ b' \mid b' \leq_f a' \}$, so $b \leq_f a'$ since $f$ is complete.
Thus, $U \leq \varphi(f)^*(a')$.

Conversely, given a presentation morphism $f : \widetilde{X} \to Y$, we construct a complete premorphism $\varphi'(f) : X \to Y$.
We define $b \leq_{\varphi'(f)} a$ if and only if $b \in f^*(a)$, and set $\varphi'(f)^\sharp_{b,a} = \pi_b \circ f^\sharp_a$.
It is complete since $f^*(a)$ is a $\cover$-ideal.
Conditions \eqref{ssm:comp}, \eqref{ssm:nat}, \eqref{ssm:top}, and \eqref{ssm:meet} of \rdefn{ssm} follow from the corresponding conditions of \rdefn{pres-map}.
To verify condition~\eqref{ssm:loc}, let $b \in f^*(a)$ and $a' \leq a$, and let $\mathcal{O}^Y_{a',a}$ be a localization at $x$.
By \axref{SE} for $X$, there exists $b' \leq b$ such that $\mathcal{O}^X_{b',b}$ is a localization at $\pi_b(f^\sharp_a(x))$.
Then the restriction of $f^\sharp_a(x)$ to $\overline{b'}$ is invertible, so \rdefn{pres-map}.\eqref{pres:loc} implies that $b' \in f^*(a')$.
Thus, $b \in L_{\varphi'(f)}(a,a')$.

\begin{lem}[phi-bij]
The constructions $\varphi$ and $\varphi'$ are mutually inverse.
\end{lem}
\begin{proof}
This follows directly from the definitions.
\end{proof}

Thus, we obtain a bijection $\psi \circ \varphi$ from the set of morphisms of saturated scheme sites $X \to Y$ to the set of morphisms of schemes $\widetilde{X} \to \widetilde{Y}$.
To show that this bijection preserves composition, it is convenient to introduce auxiliary composition functions of presentation morphisms with premorphisms of scheme sites on the right, and with morphisms of locally ringed locales on the left.

Let $f : (L,\mathcal{O}^L) \to (M,\mathcal{O}^M)$ be a morphism of locally ringed locales, and let $g : (M,\mathcal{O}^M) \to Y$ be a presentation morphism.
Define $g \circ_l f : (L,\mathcal{O}^L) \to Y$ as follows:
\begin{align*}
(g \circ_l f)^*(a) & = f^*(g^*(a)), \\
(g \circ_l f)^\sharp_a & = f^\sharp_{g^*(a)} \circ g^\sharp_a.
\end{align*}

Let $f : (L,\mathcal{O}^L) \to X$ be a presentation morphism, and let $g : X \to Y$ be a premorphism of saturated scheme sites.
We define $g \circ_r f : (L,\mathcal{O}^L) \to Y$.
For every $a \in O_L$, we set
$(g \circ_r f)^*(a) = \bigvee \{ f^*(b) \mid b \leq_g a \}$.
The homomorphism $(g \circ_r f)^\sharp_a$ is defined as the unique map that fits in the following diagram for every $b \leq_g a$:
\[ \xymatrix{ \mathcal{O}^Y_a \ar[rr]^-{(g \circ_r f)^\sharp a} \ar[d]_{g^\sharp_{b,a}} & & \mathcal{O}^L_{(g \circ_r f)^*(a)} \ar[d]^{\mathcal{O}^L_{f^*(b),(g \circ_r f)^*(a)}} \\
              \mathcal{O}^X_b \ar[rr]_-{f^\sharp_b}                                     & & \mathcal{O}^L_{f^*(b)}
            } \]

\begin{lem}[func-comp]
The operations $\circ_l$ and $\circ_r$ are well-defined and satisfy the following compatibility conditions:
\begin{itemize}
\item For every presentation morphism $f : (L,\mathcal{O}^L) \to X$ and all premorphisms $g,g' : X \to Y$, if $g \sim g'$, then $g \circ_r f = g' \circ_r f$.
\item For every complete premorphism $f : X \to Y$ and every premorphism $g : Y \to Z$ of saturated scheme sites, we have $\varphi(\overline{g \circ f}) = g \circ_r \varphi(f)$.
\item For all presentation morphisms $f : (L,\mathcal{O}^L) \to X$ and $g : \widetilde{X} \to Y$, we have $g \circ_l \psi(f) = \varphi'(g) \circ_r f$.
\item For every morphism of locally ringed locales $f : (L,\mathcal{O}^L) \to (M,\mathcal{O}^M)$ and every presentation morphism $g : (M,\mathcal{O}^M) \to Y$, we have $\psi(g \circ_l f) = \psi(g) \circ f$.
\end{itemize}
\end{lem}
\begin{proof}
This follows by a straightforward verification from the definitions.
To prove that $g \circ_r f$ is a presentation morphism, we use \rlem{loc-forall} and \rlem{map-cover} to verify condition~\eqref{pres:loc}.
The first two compatibility conditions follow from the fact that $f^*$ preserves covers (see the proof of the previous lemma).
\end{proof}

Now, we can prove that the category of saturated scheme sites embeds fully faithfully into the category of schemes:

\begin{prop}[fff]
The following mappings determine a fully faithful functor from $\cat{SatSchSite}$ to $\cat{Sch}$:
\begin{align*}
X & \mapsto \widetilde{X}, \\
f & \mapsto \psi(\varphi(f)).
\end{align*}
\end{prop}
\begin{proof}
Since $b \leq_{\overline{\mathrm{id}_X}} a$ if and only if $b \cover \{ a \}$, that is, $b \in \overline{a}$, it is straightforward to see that $\varphi(\overline{\mathrm{id}_X}) = \psi'(\mathrm{id}_{\widetilde{X}})$, which implies that $\psi \circ \varphi$ preserves identity morphisms.
The fact that it also preserves composition follows from \rlem{func-comp}:
\begin{align*}
& \psi(\varphi(\overline{g \circ f})) = \psi(g \circ_r \varphi(f)) = \psi(\varphi'(\varphi(g)) \circ_r \varphi(f)) = \\
& \psi(\varphi(g) \circ_l \psi(\varphi(f))) = \psi(\varphi(g)) \circ \psi(\varphi(f)).
\end{align*}
Thus, we indeed get a functor $\cat{SatSchSite} \to \cat{Sch}$.
\rlem{psi-bij} and \rlem{phi-bij} imply that it is fully faithful.
\end{proof}

\subsection{From schemes to scheme sites}
\label{sec:scheme-to-scheme-site}

We now construct a saturated scheme site associated to a scheme. In fact, the construction is more general and applies to any ringed locale.

First, we prove two simple lemmas about principal opens.

\begin{lem}[po-restrict]
Let $(L,\mathcal{O}^L)$ be a ringed locale.
Then for all $b,a \in O_L$ such that $b \leq a$, and every $x \in \mathcal{O}^L_a$, we have $P(a,x) \wedge b = P(b, \mathcal{O}^L_{b,a}(x))$.
\end{lem}
\begin{proof}
Recall that
\[ P(a,x) = \bigvee \{ c \in O_L \mid c \leq a, \mathcal{O}^L_{c,a}(x) \text{ is invertible} \}. \]

We first prove that
$P(b, \mathcal{O}^L_{b,a}(x)) \leq P(a,x) \wedge b$.
It suffices to show that for every $c \leq b$ such that $\mathcal{O}^L_{c,b}(\mathcal{O}^L_{b,a}(x))$ is invertible, we have $c \leq P(a,x) \wedge b$.
The inequality $c \leq b$ holds by assumption, and $c \leq P(a,x)$ follows from the fact that $\mathcal{O}^L_{c,a}(x)$ is invertible.
Thus $c \leq P(a,x) \wedge b$.

We now prove the opposite inequality.
It suffices to show that for every $c \leq a$ such that $\mathcal{O}^L_{c,a}(x)$ is invertible, we have $c \wedge b \leq P(b, \mathcal{O}^L_{b,a}(x))$.
This follows from the fact that $\mathcal{O}^L_{c \wedge b, b}(\mathcal{O}^L_{b,a}(x))$ is invertible.
\end{proof}

\begin{lem}[po-meet]
Let $(L,\mathcal{O}^L)$ be a ringed locale.
Then for all $a \in O_L$ and $x,y \in \mathcal{O}^L_a$, we have $P(a,x) \wedge P(a,y) = P(a, x y)$.
\end{lem}
\begin{proof}
This follows directly from the definition of principal opens.
\end{proof}

Let $(L,\mathcal{O}^L)$ be a ringed locale, and let $C \subseteq O_L$ be a set of affine opens such that
for all $a \in C$ and $x \in \mathcal{O}^L_a$, the principal open $P(a,x)$ also belongs to $C$.
Then we construct a saturated scheme site $S(L,\mathcal{O}^L,C)$.
Its underlying set is $C$.
The preorder is defined by declaring that
\[
b \leq a \text{ if and only if } \exists x \in \mathcal{O}^L_a, b = P(a,x).
\]
The structure presheaf is given by the restriction of $\mathcal{O}^L$ to $C$.
Condition~\axref{SE} follows from \rdefn{affine-open}.\eqref{affine:loc}.
Condition~\axref{SC} follows from \rlem{po-restrict}.

We now show that affine opens are stable under taking principal opens.

\begin{lem}[po-affine]
Let $(L,\mathcal{O}^L)$ be a ringed locale.
Then for every affine open $a \in \mathcal{O}_L$ and every $z \in \mathcal{O}^L_a$, the principal open $P(a,z)$ is also affine.
\end{lem}
\begin{proof}
We verify conditions of \rdefn{affine-open}.

\eqref{affine:loc} This follows from \rlem{loc-factor}.

\eqref{affine:surj} Let $b$ be an open such that $b \leq P(a,z)$. Then we have 
\[ b \leq \bigvee \{ c \mid c \leq b, \exists x \in \mathcal{O}^L_a, c = P(a,x) \}. \]
Thus, it is enough to show that for every $x \in \mathcal{O}^L_a$ such that $P(a,x) \leq b$, we have $P(P(a,z), \mathcal{O}^L_{P(a,z),a}(x)) \leq b$.
By \rlem{po-restrict}, we have $P(P(a,z), \mathcal{O}^L_{P(a,z),a}(x)) = P(a,x) \wedge P(a,z)$.
This implies the required inequality.

\eqref{affine:cover}
Let $x \in \mathcal{O}^L_{P(a,z)}$ and $U \subseteq O_L$ be such that $P(P(a,z),x) \leq \bigvee U$.
We need to show that the ideal of $\mathcal{O}^L_{P(P(a,z),x)}$ generated by the following set is the unit ideal:
\[ \{ y \mid \exists c \in U, \exists b, b \leq P(P(a,z),x), b \leq c, \mathcal{O}^L_{b,P(P(a,z),x)} \text{ is a localization at } y \}. \]

We can assume that $x = \mathcal{O}^L_{P(a,z),a}(x')$ for some $x' \in \mathcal{O}^L_a$.
By \rlem{po-restrict} and \rlem{po-meet}, we have $P(P(a,z),x) = P(a,x') \wedge P(a,z) = P(a, x' z)$.
By condition~\eqref{affine:cover} for $a$, the ideal of $\mathcal{O}^L_{P(a, x' z)}$ generated by the following set is the unit ideal:
\[ \{ y \mid \exists c \in U, \exists b, b \leq P(a, x' z), b \leq c, \mathcal{O}^L_{b,P(a, x' z)} \text{ is a localization at } y \}. \]

This concludes the proof since the rings $\mathcal{O}^L_{P(P(a,z),x)}$ and $\mathcal{O}^L_{P(a, x' z)}$, and the ideals described above coincide.
\end{proof}

This lemma implies that we can take $C$ to be the set of all affine opens in the construction of $S(L,\mathcal{O}^L,C)$.
In this case, we denote this scheme site by $S(L,\mathcal{O}^L)$.

\begin{prop}[surj-obj]
Every scheme $(L,\mathcal{O}^L)$ is isomorphic to $\widetilde{S(L,\mathcal{O}^L)}$.
\end{prop}
\begin{proof}
We first construct a presentation morphism $f : (L,\mathcal{O}^L) \to S(L,\mathcal{O}^L)$.
Since the set of opens of $S(L,\mathcal{O}^L)$ is a subset of $O_L$, we can define $f^*$ as the inclusion of this subset.
The homomorphism $f^\sharp_a$ is defined as the identity map.
Now, we verify conditions of \rdefn{pres-map}.

\eqref{pres:top}
Since $(L,\mathcal{O}^L)$ is a scheme, we have that $\top$ in $O_L$ is covered by the set of affine opens.

\eqref{pres:meet}
Let $a,b \in O_L$ be affine opens. We need to show that
\[ a \wedge b \leq \bigvee \{ c \in O_L \mid c \text{ is affine}, c \leq a, c \leq b \}. \]
By \rdefn{affine-open}.\eqref{affine:surj}, we know that
\[ a \wedge b \leq \bigvee \{ a' \wedge b' \mid a',b' \leq a \wedge b, \exists x \in \mathcal{O}^L_a, \exists y \in \mathcal{O}^L_b, a' = P(a,x), b' = P(b,y) \}. \]
Thus, it is enough to show that if $P(a,x),P(b,y) \leq a \wedge b$, then $P(a,x) \wedge P(b,y)$ is affine.
By \rlem{po-affine}, $P(b,y)$ is affine.
By \rlem{po-restrict}, we have $P(a,x) \wedge P(b,y) = P(P(b,y),\mathcal{O}^L_{P(b,y),a}(x))$.
Thus, by \rlem{po-affine}, $P(a,x) \wedge P(b,y)$ is also affine.

\eqref{pres:loc}
This follows immediately from the definition of principal opens.

Next, we prove that $\psi(f) : (L,\mathcal{O}^L) \to \widetilde{S(L,\mathcal{O}^L)}$ is an isomorphism.
To do this, it is enough to show that both $\psi(f)^*$ and $\psi(f)^\sharp$ are bijective.

To prove injectivity of $\psi(f)^*$, it is enough to show that $f^*$ reflects covers.
That is, we need to show that given an affine open $a \in O_L$ and a set of affine opens $U \subseteq O_L$, if $a \leq \bigvee U$, then $a$ is covered by $U$ in $S(L,\mathcal{O}^L)$.
This follows from \rdefn{affine-open}.\eqref{affine:cover} applied to $x = 1$.

To prove surjectivity of $\psi(f)^*$, we need to show that every open $b \in O_L$ is a join of affine opens.
Since $(L,\mathcal{O}^L)$ is a scheme, $\top$ is covered by a set of affine opens $\{ a_i \}_{i \in I}$.
It follows that $b = \bigvee_{i \in I} b \wedge a_i$.
Thus, it is enough to show that $b \wedge a_i$ is a join of affine opens, but this follows from \rdefn{affine-open}.\eqref{affine:surj}. 

Finally, $\psi(f)^\sharp$ is the identity map on basic opens.
It follows that it is an isomorphism on all opens by the sheaf condition.
\end{proof}

Finally, we are ready to prove the main result of this section.

\begin{thm}[main]
The categories $\cat{SatSchSite}$ and $\cat{Sch}$ are equivalent.
\end{thm}
\begin{proof}
By \rprop{fff}, we have a fully faithful functor $\cat{SatSchSite} \to \cat{Sch}$.
By \rprop{surj-obj}, this functor is essentially surjective on objects.
Thus, it is an equivalence of categories.
\end{proof}

\bibliographystyle{amsplain}
\bibliography{ref}

\end{document}